\documentclass[12pt,a4paper]{article}

\usepackage{latexsym}
\usepackage{amsfonts,amsmath}
\usepackage{graphicx}

\author{K\'aroly J. B\"or\"oczky\footnote{Supported by OTKA 75016}}
\title{Stronger versions of the Orlicz-Petty projection inequality}

\newcommand{\proof}{{\it Proof: }}

\newcommand{\R}{\mathbb{R}}

\newcommand{\theend}{\mbox{}\hfill\raisebox{-1ex}{Q.E.D.}\\}

\begin{document}

\maketitle

\newtheorem{theo}{Theorem}
\newtheorem{coro}[theo]{Corollary}
\newtheorem{lemma}[theo]{Lemma}
\newtheorem{remark}[theo]{Remark}
\newtheorem{prop}[theo]{Proposition}
\newtheorem{conj}[theo]{Conjecture}
\newtheorem{example}[theo]{Example}

\begin{abstract}
We verify a conjecture of Lutwak, Yang, Zhang about
the equality case in the Orlicz-Petty projection inequality, and  provide an
essentially optimal stability version.
\end{abstract}

The Petty projection inequality (Theorem~\ref{Petty}),  its $L_p$ extension, and its analytic counterparts,
the Zhang-Sobolev inequality \cite{Zha99} and its $L_p$ extension by
A. Cianchi, E. Lutwak, D. Yang, G. Zhang \cite{CLYZ09,LYZ02}, are fundamental affine isoperimetric and affine analytic inequalities (see in addition, e.g.,
D. Alonso-Gutierrez, J. Bastero, J. Bernu\'es \cite{ABB},
R.J.  Gardner, G. Zhang \cite{GaZ98},
C. Haberl, F.E. Schuster \cite{HaS09a,HaS09b},  C. Haberl, F.E. Schuster, J. Xio \cite{HSX},
E. Lutwak, D. Yang, G. Zhang \cite{LYZ00,LYZ06,LYZ10}, M. Ludwig \cite{Lud05,Lud10},
M. Schmuckenschl\"ager \cite{Sch95}, F.E. Schuster, T. Wannnerer \cite{ScW12},
J. Xiao \cite{Xia07}). The notion of projection body and its $L_p$ extension have found their natural context in the work of
 E. Lutwak, D. Yang, G. Zhang \cite{LYZ10}, where the authors introduced the concept of Orlicz projection body. The fundamental result of  \cite{LYZ10} is
the Orlicz-Petty projection inequality. The goal of this paper is to strengthen this latter inequality
 extending the method of E. Lutwak, D. Yang, G. Zhang \cite{LYZ10}
based on Steiner symmetrization.

When the equality case of a geometric inequality is characterized, it is a natural question how close a convex body $K$  is to the extremals if almost inequality holds for $K$ in the inequality. Precise answers to these questions are called stability versions of the original inequalities. Stability results for geometric estimates have important applications, see for example
B. Fleury, O. Gu\'edon, G.  Paouris \cite{FGP07} for the central limit theorem on convex bodies, and
D. Hug, R. Schneider \cite{HuS10} for the shape of typical cells in a Poisson hyperplane process.

Stability versions of sharp geometric inequalities have been
around since the days of Minkowski, see the survey paper by H. Groemer \cite{Gro93}
about developments until the early 1990s. Recently essentially optimal results
were obtained by N. Fusco, F. Maggi, A. Pratelli \cite{FMP08}
concerning the isoperimetric inequality, and by  A. Figalli, F. Maggi, A. Pratelli \cite{FMP09} and \cite{FMP10} for the Brunn-Minkowski inequality, see F. Maggi \cite{Mag08} for
a survey of their methods. In these papers, stability is understood in terms of volume difference
of normalised convex  bodies. In this paper we follow J. Bourgain and J. Lindenstrauss \cite{BoL88},
who used the so called Banach-Mazur distance for their result (\ref{BL}) about  projection bodies quoted below.

We write $o$ to denote the origin in $\R^n$, $u\cdot v$ to denote the scalar product of
the vectors $u$ and $v$,
$\mathcal{H}$ to denote
the $(n-1)$-dimensional Hausdorff measure,
and $[X_1,\ldots,X_k]$ to denote the convex hull
of the sets $X_1,\ldots,X_k$ in $\R^n$.
 For a non-zero $u$ in $\R^n$, let $u^\bot$ be
the orthogonal linear $(n-1)$-subspace, and let
$\pi_u$ denote the
orthogonal projection onto $u^\bot$.
In addition, let $B^n$ be
the Euclidean unit ball, and let $\kappa_n$ be its volume.
For $x\in \R^n$, $\|x\|$ denotes the Euclidean norm.
We write $A\Delta B$ to denote the symmetric difference
of the sets $A$ and $B$.

Throughout this article, a convex body in $\R^n$ is a compact convex set with non-empty interior.
 In addition, we write $\mathcal{K}_o^n$ to denote the set of convex bodies in $\R^n$
 that contain the origin in their interiors.
For a convex body $K$ in $\R^n$,
 let $h_K(u)=\max_{x\in K}x\cdot u$
denote the support function of $K$ at $u\in\R^n$,
and let $K^*$ be  be the polar of $K$, defined by
$$
K^*=\{u\in\R^n:\,h_{K}(u)\leq 1\}.
$$
Let $S_K$ be the surface area measure of $K$
on $S^{n-1}$. That is, if $\sigma$ is an open subset of $S^{n-1}$,
then $S_K(\sigma)$ is the $(n-1)$-dimensional Hausdorff-measure
of all $x\in\partial K$, where there exists an exterior unit normal lying in $\sigma$.
 Minkowski's projection body $\Pi K$ is the
 $o$-symmetric convex body whose support function is
$$
h_{\Pi K}(x)=\|x\| \cdot\mathcal{H}(\pi_x K)
=\frac12\int_{S^{n-1}}|x\cdot w|\,d S_K(w)
\mbox{ \ for $x\in\R^n\backslash o$}.
$$
We write $\Pi^*K$ to denote the polar of $\Pi K$, and
note that $V(\Pi^* K)V(K)^{n-1}$
is invariant under affine transformations of $\R^n$
(see E. Lutwak \cite{Lut93}).
Petty's projection inequality can now be stated as follows.

\begin{theo}[Petty]
\label{Petty}
If $K$ is a convex body in $\R^n$, then
$$
V(\Pi^* K)V(K)^{n-1}\leq (\kappa_n/\kappa_{n-1})^n,
$$
with equality if and only if $K$ is an ellipsoid.
\end{theo}

To define the Orlicz projection body introduced by
E. Lutwak, D. Yang, G. Zhang \cite{LYZ10},
we write $\mathcal{C}$ to denote the set of convex functions $\varphi:\,\R\to[0,\infty)$ such that $\varphi(0)=0$, and $\varphi(-t)+\varphi(t)>0$ for $t\neq 0$. In particular,
\begin{equation}
\label{phimonotone}
\mbox{every $\varphi\in\mathcal{C}$ is }
\left\{
\begin{array}{l}
\mbox{either strictly montone decreasig on $(-\infty,0]$,}\\
\mbox{or  strictly montone increasig on $[0,\infty)$.}
\end{array}\right.
\end{equation}

Let $\varphi\in\mathcal{C}$, and let $K\in\mathcal{K}^n_o$. The corresponding Orlicz projection body
$\Pi_\varphi K$ is defined in \cite{LYZ10} {\it via} its support function such that
if $x\in\R^n$, then
\begin{equation}
\label{Oprojdef}
h_{\Pi_\varphi K}(x)=\min\left\{\lambda>0:\;
\int_{S^{n-1}}\varphi\left(\frac{x\cdot w}{\lambda h_K(w)}\right)h_K(w)\,d S_K(w)
\leq nV(K)\right\}.
\end{equation}
Since the surface area measure of every open hemisphere is positive,
(\ref{phimonotone}) yields that the minimum  in (\ref{Oprojdef}) is attained at a unique $\lambda>0$.

An important special case is when $\varphi(t)=|t|^p$ for some $p\geq 1$. Then
$\Pi_\varphi K$ is the $L_p$ projection body $\Pi_p K$ introduced by
E. Lutwak, D. Yang, G. Zhang \cite{LYZ00} (using a different normalization):
\begin{equation}
\label{pdef}
h_{\Pi_p K}(x)^p=\frac1{nV(K)}
\int_{S^{n-1}}|x\cdot w|^p h_K(w)^{1-p}\,d S_K(w).
\end{equation}
In particular, if $p=1$, then
$$
\Pi_1(K)=\frac2{nV(K)} \cdot \Pi K.
$$
In addition, if $p$ tends to infinity, then we may define
the $L_\infty$ polar projection body $\Pi^*_\infty$ to be $K\cap (-K)$.

Unlike $\Pi K$, the Orlicz projection body
$\Pi_\varphi K$ is not translation invariant for a general $\varphi\in\mathcal{C}$,
and may not be $o$-symmetric. However E. Lutwak, D. Yang, G. Zhang \cite{LYZ10}
show that
\begin{equation}
\label{linearmap}
\mbox{$\Pi^*_\varphi AK=A\Pi^*_\varphi K$ holds for any $A\in{\rm GL}(n)$,
$K\in\mathcal{K}^n_o$ and $\varphi\in\mathcal{C}$}.
\end{equation}
 The following
Orlicz-Petty projection inequality is the main result of \cite{LYZ10}.

\begin{theo}[Lutwak,Yang,Zhang]
\label{PettyOrlicz}
Let $\varphi\in\mathcal{C}$.
If $K\in\mathcal{K}^n_o$, then the volume ratio
$$
\frac{V(\Pi_\varphi^* K)}{V(K)}
$$
is maximized when $K$ is an $o$-symmetric ellipsoid. If
$\varphi$ is strictly convex, then the $o$-symmetric ellipsoids
are the only maximizers.
\end{theo}

If $\varphi(t)=|t|$, which is the case of the normalized classical projection body,
 then every ellipsoid is a maximizer  in the Orlicz-Petty projection inequality 
(see Theorem~\ref{Petty}). Thus to summarize what to expect
for an arbitrary $\varphi\in\mathcal{C}$,
 E. Lutwak, D. Yang, G. Zhang \cite{LYZ10} conjecture that  every maximizer is an ellipsoid.
Here we confirm this conjecture.

\begin{theo}
\label{OrliczPettyequa}
Let $\varphi\in\mathcal{C}$.
If $K\in\mathcal{K}^n_o$ maximizes  the volume ratio $V(\Pi_\varphi^* K)/V(K)$,
then $K$ is an ellipsoid.
\end{theo}

A natural tool for stability results of affine invariant inequalities is the  Banach-Mazur distance $\delta_{\rm BM}(K,M)$
of the convex bodies $K$ and $M$ defined by
\begin{multline*}
\delta_{\rm BM}(K,M)=\min\{\lambda\geq 0:\,
K-x\subset \Phi(M-y)\subset e^\lambda(K-x)\\
\mbox{ for \ } \Phi\in{\rm GL}(n),\,x,y\in\R^n\}.
\end{multline*}
In particular, if $K$ and $M$ are $o$-symmetric, then
$x=y=o$ can be assumed. In addition, for a line $l$ passing through the origin $o$, we write $\mathcal{K}_l$ to denote
the set of $o$-symmetric convex bodies 
with axial rotational symmetry around the line $l$. 
 If $K\in \mathcal{K}_l$, then
$$
\delta_{\rm BM}(K,B^n)=\min\{\lambda\geq 0:\,E\subset K \subset e^\lambda E, \mbox{ where $E\in \mathcal{K}_l$
is an ellipsoid}\}.
$$
It follows for example
from a theorem of F. John \cite{Joh37}
that  $\delta_{\rm BM}(K,B^n)\leq \ln n$
for any convex body $K$ in $\R^n$.

We strengthen Theorem~\ref{OrliczPettyequa} as follows, where
 we set $\tilde{\varphi}(t)=\varphi(-t)+\varphi(t)$ for $\varphi\in\mathcal{C}$.

\begin{theo}
\label{OrliczPettystab}
If $\varphi\in\mathcal{C}$ and $K\in\mathcal{K}^n_o$ with $\delta=\delta_{BM}(K,B^n)$, then
$$
\frac{V(\Pi_\varphi^* K)}{V(K)}\leq
\left(1-\gamma\cdot\delta^{cn}\cdot \tilde{\varphi}(\delta^{c})\right)\cdot
\frac{V(\Pi_\varphi^* B^n)}{V(B^n)},
$$
where $c=840$ and  $\gamma>0$ depends on $n$ and $\varphi$.
\end{theo}

Next we discuss what Theorem~\ref{OrliczPettystab} yields for  Petty's projection inequality.

\begin{coro}
\label{Pettystab}
If $K$ is a convex body in $\R^n$ with $\delta=\delta_{BM}(K,B^n)$, then
$$
V(\Pi^* K)V(K)^{n-1}\leq \left(1-\gamma\cdot\delta^{cn}\right)
(\kappa_n/\kappa_{n-1})^n
$$
where $c=1680$ and $\gamma>0$ depends only on $n$.
\end{coro}

The example below shows that
the exponent $cn$ for an absolute constant $c>0$ is
of optimal order.
G. Ambrus and the author \cite{AmB} recently proved Corollary~\ref{Pettystab}
with an exponent of the form $cn^3$ instead of the optimal $cn$. \\

\noindent{\bf Example }
Let $K=[B^n,\pm(1+\varepsilon)v)]$ for
some $v\in S^{n-1}$.
In this case, the Banach-Mazur distance of $K$
from any ellipsoid is at least
$\varepsilon/2$,
and
$$
V(\Pi^* K)V(K)^{n-1}\geq (1-\gamma_0\varepsilon^{\frac{n+1}2})(\kappa_n/\kappa_{n-1})^n,
$$
where $\gamma_0>0$ depends only on $n$. \\

As a related result, J. Bourgain and J. Lindenstrauss \cite{BoL88}
proved that if $K$ and $M$ are $o$-symmetric convex bodies in $\R^n$, then
\begin{equation}
\label{BL}
\delta_{BM}(\Pi K,\Pi M)\geq \gamma\cdot \delta_{BM}(K,M)^{n(n+5)/2}
\end{equation}
where $\gamma>0$ depends only on $n$, and they
conjectured that the optimal order of the exponent
is $cn$ for an absolute constant $c>0$. 
The exponent in (\ref{BL}) has been slightly improved by S. Campi \cite{Cam88} if $n=3$, and by
M. Kiderlen \cite{Kid08} for any $n$, but the conjecture is still wide open.
Corollary~\ref{Pettystab}
is in accordance with this conjecture of J. Bourgain and J. Lindenstrauss in the case when $M$ is an ellipsoid.
Actually, if $K$ and $M$ are not $o$-symmetric then
their projection bodies may coincide
even if $\delta_{BM}(K,M)\neq 0$ (see R. Schneider \cite{Sch67}).

If $\varphi$ is strictly convex, then
 E. Lutwak, D. Yang, G. Zhang \cite{LYZ10}  proved that
the $o$-symmetric ellipsoids are the only maximizers
 in the Orlicz-Petty projection inequality (see Theorem~\ref{PettyOrlicz}).
We prove a stability version of this statement for even $\varphi$. For $K\in\mathcal{K}_o^n$, let
$$
\delta^*_{\rm EL}(K)=\min\{\lambda\geq 0:\,
E\subset K\subset e^\lambda E\mbox{ for some $o$-symmetric ellipsoid $E$}\}.
$$
Since $\delta^*_{\rm EL}(K)$ becomes arbitrary large if $K$ is translated in a way such that the origin gets close to $\partial K$, it is more natural to consider
$$
\delta_{\rm EL}(K)=\min\{1,\delta^*_{\rm EL}(K)\}.
$$

\begin{theo}
\label{OrliczPettystrictstab}
Let  $\varphi\in\mathcal{C}$ be even such that $\varphi''(t)$ is continuous
and positive for $t>0$.
If $K\in\mathcal{K}^n_o$ with $\delta=\delta_{\rm EL}(K)$, then
$$
\frac{V(\Pi^*_\varphi K)}{V(K)}\leq\left(1-\gamma\cdot
\delta^{cn}\cdot\varphi(\delta^{c})\right)
\cdot\frac{V(\Pi^*_\varphi B^n)}{V(B^n)},
$$
where  $c=2520$ and $\gamma>0$ depends only on $n$ and $\varphi$.
\end{theo}

Under the conditions of Theorem~\ref{OrliczPettystrictstab}, let $K\in\mathcal{K}^n_o$ be such that
$V(\Pi^*_\varphi K)/V(K)$ is very close to $V(\Pi^*_\varphi B^n)/V(B^n)$. Then
 Theorem~\ref{OrliczPettystab} yields that there exists a translate $K'$ of $K$ such that $\delta_{\rm EL}(K')$
is small, while  Theorem~\ref{OrliczPettystrictstab} implies that already $\delta_{\rm EL}(K)$
is small.

For the $L_p$ projection body for $p>1$, and for $c=2520$, we have
$$
\frac{V(\Pi^*_p K)}{V(K)}\leq\left(1-\gamma\cdot
\delta_{\rm EL}(K)^{c(n+p)}\right)
\cdot\frac{V(\Pi^*_p B^n)}{V(B^n)}.
$$
Here the order of the error term
gets smaller and smaller as $p$ grows. It is not surprising, because
$\Pi_\infty^*(K)=K\cap (-K)$
for $K\in\mathcal{K}_o^n$, and hence $V(\Pi_\infty^* K)/V(K)$ is maximized
by any $o$-symmetric convex body $K$.

Our arguments to prove Theorems~\ref{OrliczPettyequa}, \ref{OrliczPettystab}
and \ref{OrliczPettystrictstab} are based on Steiner symmetrization, and are variations of the method developed in  E. Lutwak, D. Yang, G. Zhang  \cite{LYZ10}.
The novel ideas to prove Theorems~\ref{OrliczPettyequa} and \ref{OrliczPettystab}
are to compare shadow boundaries in two suitable independent directions, and
to reduce the problem to convex bodies with axial rotational symmetry around $\R u$
for a $u\in S^{n-1}$. In the latter case, the shadow boundaries parallel to $u$
and orthogonal to $u$ are well understood, which makes it possible to
perform explicit caculations.

For Theorem~\ref{OrliczPettystab},  the proof of the reduction
to convex bodies with axial rotational symmetry is rather technical,
so the argument for the corresponding statement
Theorem~\ref{equatorpole} is deferred to Section~\ref{secequatorpole}.

We note that W. Blaschke \cite{Bla16} characterized ellipsoids as the only convex bodies
such that every shadow boundary is contained in some hyperplane. A stability version of this statement was proved by P.M. Gruber \cite{Gru97}.

\section{Some facts about convex bodies}
\label{secclassical}

Unless we provide specific references,
 the results reviewed in this section are discussed in
 the monographs by
T. Bonnesen, W. Fenchel \cite{BoF87},
P.M. Gruber \cite{Gru07}, and R. Schneider \cite{Sch93}.
We note that the $L_\infty$-metric on the restriction of the support functions
to $S^{n-1}$ endows the space of convex bodies with the so-called Hausdorff metric.
It is well-known that volume is continuous with respect to this metric, and
Lemma~2.3 in E. Lutwak, D. Yang, G. Zhang \cite{LYZ10} says that the polar Orlicz projection body is also continuous for fixed $\varphi\in\mathcal{C}$.

We say that a convex body $M$  in $\R^n$, $n\geq 3$,
is smooth if the tangent hyperplane is unique
at every boundary point, and we say that $M$ is strictly convex if
every tangent hyperplane intersects $M$ only
in one point.

Let $K$ be a convex body in $\R^n$.
For $v\in S^{n-1}$, let $S_vK$
denote the Steiner symmetral of $K$ with respect to $v^\bot$.
In particular, if $f,g$ are the concave real functions on $\pi_vK$ such
that
$$
K=\{y+tv:y\in \pi_vK,-g(y)\leq t\leq f(y)\},
$$
then
\begin{equation}
\label{Steinerdef}
S_vK=\{y+tv:y\in \pi_vK,\;|t|\leq \mbox{$\frac{f(y)+g(y)}2\,$}\}.
\end{equation}
Fubini's theorem yields that $V(S_vK)=V(K)$.
It is known that for any convex body $K$, there is
a sequence of Steiner symmetrizations whose
limit is a ball (of volume $V(K)$).

Next there exists
a sequence of Steiner symmetrizations
 with respect to $(n-1)$-subspaces containing the line $\R v$
 such that their limit is a convex body $R_vK$
 whose axis of rotational symmetry is $\R v$.
 This $R_vK$ is the Schwarz rounding of $K$
 with respect to $v$. In particular  a hyperplane $H$ intersects ${\rm int}\,K$
if and only if it intersects ${\rm int}\,R_vK$, and $\mathcal{H}(H\cap K)=\mathcal{H}(H\cap R_vK)$ in this case.

For our arguments, it is crucial to have
a basic understanding of the boundaries
of convex bodies.
For $x\in \partial K$, let $w_x$ be a unit exterior normal
to $\partial K$ at $x$. The following two
well-known properties
are consequences of the fact that Lipschitz functions
are almost everywhere differentiable.
\begin{description}
\item{(i)} $w_x$ is uniquely determined at
$\mathcal{H}$ almost all $x\in \partial K$.
\item{(ii)} The supporting hyperplane
with exterior normal vector $u$ intersects
$\partial K$ in a unique point
for almost all $u\in S^{n-1}$.
\end{description}

The shadow boundary $\Xi_{u,K}$ of $K$ with respect to
a $u\in\R^n\backslash o$
is the family of all $x\in\partial K$ such that
the line $x+\R u$ is tangent to $K$.
In addition  we
call the shadow boundary $\Xi_{u,K}$ thin
if it contains no segment parallel to $u$.
According to
G. Ewald, D.G. Larman, C.A. Rogers \cite{RLR70}, we have

\begin{theo}[Ewald-Larman-Rogers]
\label{ELR}
If $K$ is a convex body in $\R^n$, then
 the shadow boundary $\Xi_{u,K}$ is thin
for $\mathcal{H}$-almost all $u\in S^{n-1}$.
\end{theo}

If a connected Borel $U\subset \partial K$ is disjoint from the shadow boundary with respect a $v\in S^{n-1}$, then for any measurable 
$\psi:\pi_v(U)\to\R$, we have
\begin{equation}
\label{RadNik}
\int_{\pi_v(U)}\psi(y)\,dy=\int_{U}\psi(\pi_v x)| v\cdot w_x|\,dx.
\end{equation}

If $K\in\mathcal{K}_o^n$, then let
$\varrho_K$ be the radial function of $K$ on $S^{n-1}$, defined such that
$\varrho_K(v)\,v\in\partial M$ for $v\in S^{n-1}$. It
follows that
\begin{equation}
\label{radialvol}
V(K)=\int_{S^{n-1}}\frac{\varrho_K(w)^n}n\,d\mathcal{H}(w).
\end{equation}
In addition for the polar $K^*$ of $K$, and $v\in S^{n-1}$, we have
\begin{equation}
\label{radialh}
\varrho_{K^*}(v)=h_K(v)^{-1}.
\end{equation}

We say that a convex body $M$ is in
isotropic position, if $V(M)=1$, the centroid
of $M$ is the origin, and there exists $L_M>0$
such that
$$
\int_M (w\cdot x)^2\,dx=L_M
\mbox{ \ for any $w\in S^{n-1}$}
$$
(see A. Giannopoulos \cite{Gian}, A. Giannopoulos, V.D. Milman \cite{GiM04}
and V.D. Milman, A. Pajor \cite{MiP89} for main properties).
Any convex body $K$ has an affine image
$M$ that is
in isotropic position, and we set $L_K=L_M$. We also note that if $E$ is an $o$-symmetric ellipsoid in $\R^n$, then
for any $w\in S^{n-1}$, we have
\begin{equation}
\label{ellipsmoment}
\int_E( w\cdot x)^2\,dx=h_E(w)^2V(E)^{\frac{n+2}n}L_{B^n}.
\end{equation}

Let $\varphi\in\mathcal{C}$, and let $K\in\mathcal{K}^n_o$. We collect some additional properties
of the Orlicz projection body.
The cone volume measure $V_K$ associated to $K$ on $S^{n-1}$
defined by $d\,V_K(w)=\frac{h_K(w)}{nV(K)}\,d S_K(w)$  is a probability measure
whose study was initiated by
M. Gromov, V. Milman \cite{GrM87} (see say A. Naor \cite{Nao07} for recent applications).
The definition (\ref{Oprojdef}) of
$\Pi_\varphi K$ yields
(see Lemma~2.1 in E. Lutwak, D. Yang, G. Zhang \cite{LYZ10})  that
\begin{equation}
\label{O*projdef}
\mbox{$x\in\Pi^*_\varphi K$ if and only if }
\int_{S^{n-1}}\varphi\left(\frac{x\cdot w}{h_K(w)}\right)\,d V_K(w)\leq 1.
\end{equation}

\section{Characterizing the equality case in
the Orlicz-Petty projection inequality}
\label{secPettyOrliczequa}

Our method is an extension of the argument by
E. Lutwak, D. Yang, G. Zhang \cite{LYZ10} to prove
the Orlicz-Petty projection inequality, Theorem~\ref{PettyOrlicz}, using Steiner symmetrization.
The core of the argument of \cite{LYZ10} is Corollary 3.1,  and here we also include
a consequence of Corollary 3.1 from \cite{LYZ10} for Schwarz rounding.

\begin{lemma}[Lutwak,Yang,Zhang]
\label{lem14}
If $\varphi\in\mathcal{C}$, $K\in\mathcal{K}^n_o$
 and $v\in S^{n-1}$, then
$$
S_v\Pi_\varphi^*K\subset \Pi_\varphi^*S_vK.
$$
 In particular, $V(\Pi_\varphi^*S_vK)\geq V(\Pi_\varphi^*K)$
and $V(\Pi_\varphi^*R_vK)\geq V(\Pi_\varphi^*K)$.
\end{lemma}

We recall various facts  from \cite{LYZ10} that lead to the proof of Lemma~\ref{lem14},
because we need them in the sequel.
We note that a concave function is almost everywhere differentiable
on convex sets.

Let $\varphi\in\mathcal{C}$, let $K\in\mathcal{K}^n_o$,
 and let $v$ be a unit vector in $\R^n$.
We write $w_x$ to denote an exterior unit normal
at some $x\in\partial K$.  In addition, we frequently write an $x\in\R^n$  in
the form $x=(y,t)$ if $x=y+tv$
for $y\in v^\bot$ and $t\in\R$.
If $h$ is a concave function on $\pi_v({\rm int}K)$, then
we define
$$
\langle h\rangle(z)=h(z)-z\cdot\nabla h(z)
\mbox{ \ for  $z\in\pi_v({\rm int}K)$ where $\nabla h(z)$ exists}.
$$
If $\mu_1,\mu_2>0$, and $h_1,h_2$ are concave functions on $\pi_v({\rm int}K)$, then
$$
\langle\mu_1h_1+\mu_2 h_2\rangle=\mu_1\langle h_1\rangle
+\mu_2\langle h_2\rangle.
$$

Let $f,g$ denote the concave real functions on $\pi_vK$ such
that
$$
K=\{(y,t):y\in \pi_vK,-g(y)\leq t\leq f(y)\}.
$$
If $x=(z,f(z))\in\partial K$ and $\tilde{x}=(z,-g(z))\in\partial K$
for a $z\in\pi_v({\rm int}K)$,
and both $f$ and $g$ are differentiable at $z$, then
\begin{eqnarray}
\label{ux}
w_x&=&\left(\frac{-\nabla f(z)}{\sqrt{1+\|\nabla f(z)\|^2}},
\frac{1}{\sqrt{1+\|\nabla f(z)\|^2}}\right)\\
\label{tildeux}
w_{\tilde{x}}&=&\left(\frac{-\nabla g(z)}{\sqrt{1+\|\nabla g(z)\|^2}},
\frac{-1}{\sqrt{1+\|\nabla g(z)\|^2}}\right).
\end{eqnarray}
From this, we deduce that for any $(y,t)\in\R^n$, we have
\begin{equation}
\label{lrangle}
\begin{array}{rcl}
(y,t)\cdot w_x&=&
(- y\cdot \nabla f(z)+t)\cdot(v\cdot w_x)\\
(y,t)\cdot w_{\tilde{x}}&=&
(-y\cdot \nabla g(z)-t)\cdot (v\cdot w_{\tilde{x}})\\
h_K(w_x)&=&(z,f(z))\cdot w_x=\langle f\rangle(z)\cdot(v\cdot w_x)\\
h_K(w_{\tilde{x}}))&=&(z,-g(z))\cdot w_x=-\langle g\rangle(z)
\cdot (v\cdot w_{\tilde{x}})
\end{array}
\end{equation}
Since for any $u\in\R^n$, the definitions of the cone volume measure
and the surface area measure yield that
$$
nV(K)\int_{S^{n-1}}\varphi\left(\frac{u\cdot w}{h_K(w)}\right)\,d V_K(w)=
\int_{\partial K}\varphi\left(\frac{u\cdot w_x}{h_K(w_x)}\right)h_K(w_x)\,d \mathcal{H}(x),
$$
we deduce from (\ref{RadNik}) and  (\ref{lrangle}) the following formula, which is Lemma~3.1 in \cite{LYZ10}.
We note that Lemma~3.1 in \cite{LYZ10} assumes that $\Xi_{v,K}$ is thin, but only uses this property to ensure that the corresponding integral over $\Xi_{v,K}$ is zero.

\begin{lemma}[Lutwak,Yang,Zhang]
\label{lem12}
Using the notation as above, if $\mathcal{H}(\Xi_{v,K})=0$
and $(y,t)\in\R^n$, then
\begin{align*}
&nV(K)\int_{S^{n-1}}\varphi\left(\frac{(y,t)\cdot w}{h_K(w)}\right)\,d V_K(w)=\\
&
\int_{\pi_v K}
\varphi\left(\frac{t- y\cdot\nabla f(z)}{\langle f\rangle(z)}\right)\langle f\rangle(z)\,dz
+
\int_{\pi_v K}
\varphi\left(\frac{-t- y\cdot\nabla g(z)}{\langle g\rangle(z)}\right)\langle g\rangle(z)\,dz.
\end{align*}
\end{lemma}

We continue to use the notation of Lemma~\ref{lem12} and the
condition $\mathcal{H}(\Xi_{v,K})=0$.
If $(y,t),(y,-s)\in\partial\Pi_\varphi^*K$ for $t>-s$, then it follows from (\ref{O*projdef}) that
$$
\frac12\left(\int_{S^{n-1}}\varphi\left(\frac{(y,t)\cdot w}{h_K(w)}\right)\,d V_K(w)
+\int_{S^{n-1}}\varphi\left(\frac{(y,-s)\cdot w}{h_K(w)}\right)\,d V_K(w)
\right)=1.
$$
Therefore (\ref{Steinerdef}) and Lemma~\ref{lem12} yield that
 for $(y,\frac12(t+s))\in \partial S_v\Pi_\varphi^*K$, we have
\begin{align}
\label{long1}
&nV(K)\left[1-\int_{S^{n-1}}\varphi\left(\frac{(y,\frac12(t+s))\cdot w}{h_{S_vK}(w)}\right)\,d V_{S_vK}(w)\right]=\\
\nonumber
&\frac1{2}\int_{\pi_v K}
\varphi\left(\frac{t- y\cdot\nabla f(z)}{\langle f\rangle(z)}\right)\langle f\rangle(z)\,dz
+ \frac12\int_{\pi_v K}
\varphi\left(\frac{s- y\cdot\nabla g(z)}{\langle g\rangle(z)}\right)\langle g\rangle(z)\,dz\\
\label{long2}
&-\int_{\pi_v K}
\varphi\left(\frac{\frac{t}2+\frac{s}2- \frac{y\cdot\nabla f(z)}2-\frac{y\cdot\nabla g(z)}2}
{\frac{\langle f\rangle(z)}2+\frac{\langle g\rangle(z)}2}\right)
\left(\frac{\langle f\rangle (z)}2+\frac{\langle g\rangle(z)}2\right)\,dz\\
\label{long3}
&-\int_{\pi_v K}
\varphi\left(\frac{-\frac{t}2-\frac{s}2- \frac{y\cdot\nabla f(z)}2-\frac{y\cdot\nabla g(z)}2}
{\frac{\langle f\rangle(z)}2+\frac{\langle g\rangle(z)}2}\right)
\left(\frac{\langle f\rangle(z)}2+\frac{\langle g\rangle(z)}2\right)\,dz\\
\nonumber
&+\frac12\int_{\pi_v K}
\varphi\left(\frac{-t- y\cdot\nabla f(z)}{\langle f\rangle(z)}\right)\langle f\rangle(z)\,dz
+ \frac12\int_{\pi_v K}
\varphi\left(\frac{-s- y\cdot\nabla g(z)}{\langle g\rangle(z)}\right)\langle g\rangle(z)\,dz.
\end{align}

If $\varphi\in\mathcal{C}$, $\alpha,\beta>0$, and  $a,b\in \R$, then the convexity
of $\varphi$ yields that
\begin{equation}
\label{phiconvex}
\mbox{$\alpha\varphi\left(\frac{a}{\alpha}\right)+\beta\varphi\left(\frac{b}{\beta}\right)
\geq (\alpha+\beta)\varphi\left(\frac{a+b}{\alpha+\beta}\right).$}
\end{equation}
If in addition $a\cdot b<0$, then we deduce from $\varphi(0)=0$
and (\ref{phimonotone}) that
 \begin{equation}
\label{phiconvex0}
\mbox{$\alpha\varphi\left(\frac{a}{\alpha}\right)+\beta\varphi\left(\frac{b}{\beta}\right)
> (\alpha+\beta)\varphi\left(\frac{a+b}{\alpha+\beta}\right).$}
\end{equation}
Applying (\ref{phiconvex}) in (\ref{long2}) and (\ref{long3}) shows that
 \begin{equation}
\label{inequality}
\int_{S^{n-1}}\varphi\left(\frac{(y,\frac12(t+s))\cdot w}{h_{S_vK}(w)}\right)\,d V_{S_vK}(w)\leq 1
\end{equation}
in (\ref{long1}). We conclude
$(y,\frac12(t+s))\in \Pi_\varphi^*S_v K$ from (\ref{O*projdef}),
and in turn Lemma~\ref{lem14} in the case when $\mathcal{H}(\Xi_{v,K})=0$.

So far we have just copied the argument of E. Lutwak, D. Yang, G. Zhang \cite{LYZ10}.
We take a different route only for analyzing the equality case
in Lemma~\ref{PettySteiner},
using (\ref{phiconvex0}) instead of (\ref{phiconvex}) at an appropriate place.

For a convex body $K$ in $\R^n$ and $u\in\R^n\backslash o$,
let $\Xi^+_{u,K}$ and $\Xi^-_{u,K}$ be the set of $x\in\partial K$
where all exterior unit normals have positive and negative, respectively,
scalar product with $u$. In particular,
if $\Xi_{u,K}$ is thin, then
\begin{equation}
\label{thinshadow}
\mbox{any $x\in\Xi_{u,K}$ lies in the closures of both
$\Xi_{u,K}^+$ and $\Xi_{u,K}^-$. }
\end{equation}

\begin{lemma}
\label{PettySteiner}
Let $\varphi\in\mathcal{C}$, let $K\in\mathcal{K}^n_o$, and let
$u,\tilde{u}\in\partial \Pi_\varphi^*K$  and $v\in S^{n-1}$
such that  $u$ and $\tilde{u}$
are independent, both $\Xi_{u,K}$ and
$\Xi_{\tilde{u},K}$ are thin,   $v$ is parallel to
$u-\tilde{u}$, and $\mathcal{H}(\Xi_{v,K})=0$.
If $V(\Pi_\varphi^*\mathcal{S}_vK)= V(\Pi_\varphi^*K)$,  then
$\pi_v\Xi_{u,K}=\pi_v\Xi_{\tilde{u},K}$.
\end{lemma}
\proof
Using the notation
of (\ref{long1}) with $u=(y,t)$ and $\tilde{u}=(y,-s)$,
 we write $w(z)$ and $\tilde{w}(z)$ to denote an exterior
 unit normal vector to $\partial K$ at
$(z,f(z))$ and $(z,g(z))$, respectively,  for any $z\in\pi_v({\rm int} K)$.
Since we have equality in (\ref{inequality}),
it follows from (\ref{lrangle}), (\ref{long1}) and (\ref{phiconvex0}) that
$ (u\cdot w(z))\cdot( \tilde{u}\cdot\tilde{w}(z))\geq 0$
and $( u\cdot\tilde{w}(z))\cdot ( \tilde{u}\cdot w(z))\geq 0$
for $\mathcal{H}$-almost all $z\in\pi_v({\rm int} K)$. We conclude
by continuity that if both $(z,f(z))$ and $(z,-g(z))$  are smooth points of $\partial K$ 
for a $z\in\pi_v({\rm int} K)$, then
\begin{equation}
\label{doublesmooth}
\mbox{$ (u\cdot w(z))\cdot( \tilde{u}\cdot\tilde{w}(z))\geq 0$
and $( u\cdot\tilde{w}(z))\cdot ( \tilde{u}\cdot w(z))\geq 0$}.
\end{equation} 
If $(z,f(z))$ and $(z,-g(z))$ are both smooth points of $\partial K$ 
for a $z\in\pi_v({\rm int} K)$, then we say that they are the double smooth twins of each other.
 In particular,
$\mathcal{H}$-almost all points of $\partial K$ have a double smooth twin
by $\mathcal{H}(\Xi_{v,K})=0$.
 
It follows from by (\ref{thinshadow}) and $\mathcal{H}(\Xi_{\tilde{u},K})=0$ that for any $x\in\Xi_{u,K}$, we may choose sequences
$\{x_n\}\subset\Xi_{u,K}^+$ and
$\{y_n\}\subset\Xi_{u,K}^-$ 
tending to $x$ such that $\pi_vx_n,\pi_vy_n\not\in\pi_v\Xi_{\tilde{u},K}$, and $x_n$ and $y_n$
have double smooth twins $\tilde{x}_n$ and $\tilde{y}_n$, respectively.
Thus the sequences $\{\tilde{x}_n\}$ and $\{\tilde{y}_n\}$ tend
to the same $y\in\partial K$, which readily satisfies $\pi_vy=\pi_vx$.
We have
$\{\tilde{x}_n\}\subset\Xi_{\tilde{u},K}^+$ and
$\{\tilde{y}_n\}\subset\Xi_{\tilde{u},K}^-$ by (\ref{doublesmooth}), $\pi_v\tilde{x}_n=\pi_vx_n\not\in\pi_v\Xi_{\tilde{u},K}$
and $\pi_v\tilde{y}_n=\pi_vy_n\not\in\pi_v\Xi_{\tilde{u},K}$.
Therefore $y\in \Xi_{\tilde{u},K}$. We deduce
$\pi_v\Xi_{u,K}\subset\pi_v\Xi_{\tilde{u},K}$, and in turn
$\pi_v\Xi_{\tilde{u},K}\subset\pi_v\Xi_{u,K}$ by an analogous argument.
\theend

In our argument, we reduce the problem
to convex bodies with axial rotational symmetry.
Concerning their boundary structure, we use the following
simple observation.

\begin{lemma}
\label{axial}
If $K$ is a convex body in $\R^n$ such that the line $l$
is an axis of rotational symmetry, and the
line $l_0$ intersects $\partial K$ in a  segment,
then  either $l_0$ is parallel to $l$, or
$l_0$ intersects $l$.
\end{lemma}
\proof For any $x\in K$, we write $\varrho(x)$ to denote the radius
of the section of $K$ by the hyperplane passing through $x$
and orthogonal to $l$, where $\varrho(x)=0$
if the section is just the point $x$.

Let $l_0$ intersect $\partial K$
in the segment $[p,q]$, and let $m$ be the midpoint of $[p,q]$.
We write $p',q',m'$ to denote the
orthogonal projections of $p,q,m$ respectively,
onto $l$. It follows that
\begin{eqnarray*}
\varrho(m)&\geq &\mbox{$\frac12(\varrho(p)+\varrho(q))=
\frac12(\|p-p'\|+\|q-q'\|)$}\\
&\geq&
\mbox{$\|\frac12(p-p')+\frac12(q-q')\|=\|m-m'\|$}.
\end{eqnarray*}
Since $m\in\partial K$, we have $\varrho(m)=\|m-m'\|$,
and hence the equality case of the triangle inequality yields that
$p-p'$ and $q-q'$ are parallel.
Therefore $l$ and $l_0$ are contained in a two-dimensional
affine subspace.
\theend

\noindent{\it Proof of Theorem~\ref{OrliczPettyequa}: }
It is equivalent to show  that
we have strict inequality in the
Orlicz-Petty projection inequality if $K$ is not an ellipsoid.
Let us assume this, and that $K$ is in isotropic position.
It is sufficient to prove that there
exist a unit vector $v$, and a convex body $M$ with $V(M)=1$ such that
$$
V(\Pi_\varphi^*K)\leq V(\Pi_\varphi^*M)<V(\Pi_\varphi^*\mathcal{S}_vM).
$$
The idea is to reduce the problem to bodies with axial rotational
symmetry because in this way we will have two shadow
boundaries that are contained in some hyperplanes.

Since $K$ is not a ball of center $o$, $h_K$ is not constant,
thus we may assume that for some $p\in S^{n-1}$, we have
$$
h_K(p)^2L_{B^n}\neq L_K=\int_K (p\cdot x)^2\,dx.
$$
It follows from (ii) in Section~\ref{secclassical} that we may assume
that the supporting hyperplanes with exterior normals
$p$ and $-p$ intersect $K$ in one point.

Let  $K_1$ be the Schwarz rounding
of $K$ with respect to $\R p$. In particular
$V(K_1)=V(K)=1$, $h_{K_1}(p)=h_K(p)$ and Fubini's theorem yields
$$
\int_{K_1}( p,x)^2\,dx=\int_K( p\cdot x)^2\,dx
\neq h_{K_1}(p)^2L_{B^n}.
$$
Therefore $K_1$ is not an ellipsoid according to
(\ref{ellipsmoment}), and the supporting hyperplanes with exterior normals
$p$ and $-p$ intersect $K_1$ in one point.
In particular if $q\in S^{n-1}\cap p^\bot$,
then $\Xi_{q,K_1}=q^\bot\cap\partial K_1$ is thin.
We fix a $q\in S^{n-1}\cap p^\bot$.
Since $K_1$ is
not an ellipsoid, $\Xi_{q,K_1}$ is not the relative boundary
of some $(n-1)$-ellipsoid.\\

\noindent {\bf Case 1 } $\Xi_{p,K_1}$ is thin\\
In this case, $\Xi_{p,K_1}$ is the relative boundary
of some $(n-1)$-ball. Choose $t_1,s_1>0$ such that
$u_1=t_1p\in\partial \Pi K_1$ and $\tilde{u}_1=s_1q\in\partial \Pi K_1$,
and let $v_1=(u_1-\tilde{u}_1)$.
It follows from  Lemma~\ref{axial} that
$\Xi_{v_1,K_1}$ contains at most two segments parallel to $v_1$, and
hence its $\mathcal{H}$-measure is zero.
We have already seen that $\Xi_{\tilde{u}_1,K_1}=\Xi_{q,K_1}$
is thin, therefore we may apply Lemma~\ref{PettySteiner}
to $u_1$, $\tilde{u}_1$, $v_1$.
Since
$\pi_{v_1}\Xi_{u_1,K_1}$ is the relative boundary
of some $(n-1)$-ellipsoid, and
$\pi_{v_1}\Xi_{\tilde{u}_1,K_1}$ is not,
we deduce from  Lemma~\ref{PettySteiner} that
$$
V(\Pi^*K)\leq V(\Pi^*K_1)<V(\Pi^*\mathcal{S}_{v_1}K_1).
$$

\noindent {\bf Case 2 } $\Xi_{p,K_1}$ is not thin\\
For some $\varrho,\alpha>0$, there exists
 a segment of length $\alpha$ parallel to $p$
such that $\Xi_{p,K_1}$ is the Minkowski sum of the
segment and the
relative boundary of the $(n-1)$-ball of radius $\varrho$ centered at $o$
in $p^{\bot}$. Let $K_2$ be the Schwarz rounding of $K_1$
with respect to $\R q$, and hence $\Xi_{p,K_2}$ and $\Xi_{q,K_2}$ are both thin.

For $t\in\R$, let
$$
H(q,t)=q^\bot+tq.
$$
If $\tau\in (0,\varrho)$, then
$$
\mathcal{H}(H(q,\varrho-\tau)\cap K_2)=\mathcal{H}(H(q,\varrho-\tau)\cap K_1)>\alpha\sqrt{\varrho}\kappa_{n-2}\cdot \tau^{\frac{n-2}2}.
$$
If $K_2$ were an ellipsoid, then there would exist a $\gamma>0$ depending on $K_2$ such
$\mathcal{H}(H(q,\varrho-\tau)\cap K_2)<\gamma\cdot \tau^{\frac{n-1}2}$ for $\tau\in (0,\varrho)$, therefore $K_2$ is not an ellipsoid.
 Now we choose $t_2,s_2>0$ such that
$u_2=t_2q\in\partial \Pi K_2$ and $\tilde{u}_2=s_2p\in\partial \Pi K_2$,
and let $v_2=(u_2-\tilde{u}_2)/\|u_2-\tilde{u}_2\|$.
An argument as above using Lemma~\ref{PettySteiner} yields
$$
V(\Pi^*K)\leq V(\Pi^*K_1)
\leq V(\Pi^*K_2)<V(\Pi^*\mathcal{S}_{v_2}K_2).
\mbox{ \ } \theend
$$

\section{Proof of Theorem~\ref{OrliczPettystab}}
\label{secOrliczPettystab}

The proof  is a delicate analysis of the argument of Theorem~\ref{OrliczPettyequa}.
For example, we need a stability version of (\ref{phiconvex0}).

\begin{lemma}
\label{phiaround0}
If $\varphi\in\mathcal{C}$, $\alpha,\beta,\omega>0$, and  $a,b\in \R$ such that
$a\cdot b<0$, and $\frac{|a|}{\alpha},\frac{|b|}{\beta}\geq \omega$, then
\begin{equation}
\label{phiaround0eq}
\mbox{$\alpha\varphi\left(\frac{a}{\alpha}\right)+\beta\varphi\left(\frac{b}{\beta}\right)
-(\alpha+\beta)\varphi\left(\frac{a+b}{\alpha+\beta}\right)\geq
\frac{\min\{|a|,|b|\}}{\omega}\cdot(\varphi(-\omega)+\varphi(\omega)).$}
\end{equation}
\end{lemma}
\proof We write $\Omega$ to denote the left hand side of (\ref{phiaround0eq}).
If $\mu\geq 1$ and $t\in\R$, then the convexity of $\varphi$
and $\varphi(0)=0$ yield
\begin{equation}
\label{convexstretch}
\varphi(\mu\,t)\geq \mu\cdot\varphi(t).
\end{equation}
We may assume that $a\geq -b>0$. In particular
$0\leq  \frac{a+b}{\alpha+\beta}< \frac{a}{\alpha}$, and we deduce from
(\ref{convexstretch}) the estimate
$$
\varphi\left(\frac{a+b}{\alpha+\beta}\right)\leq
\frac{\alpha(a+b)}{a(\alpha+\beta)}\cdot\varphi\left(\frac{a}{\alpha}\right).
$$
It follows from this inequality and (\ref{convexstretch}) that
\begin{eqnarray*}
\Omega&\geq &\alpha\varphi\left(\frac{a}{\alpha}\right)+
\beta\varphi\left(\frac{b}{\beta}\right)-
\frac{\alpha(a+b)}{a}\cdot\varphi\left(\frac{a}{\alpha}\right)\\
&=&\beta\varphi\left(\frac{b}{\beta}\right)+
\frac{\alpha (-b)}{a}\cdot\varphi\left(\frac{a}{\alpha}\right)
\geq \frac{|b|}{\omega}\varphi(-\omega)+\frac{|b|}{\omega}\varphi(\omega).
\mbox{ \ }\theend.
\end{eqnarray*}

We also need the stability version Lemma~\ref{O*projdefstab} of (\ref{O*projdef}). Let $c_\varphi>0$ be defined by 
$\max\{\varphi(-c_\varphi),\varphi(c_\varphi)\}=1$ for $\varphi\in\mathcal{C}$. According to
 Lemma~2.2 by  E. Lutwak, D. Yang, G. Zhang \cite{LYZ10} stated for the Orlicz projection body, if $rB^n\subset K\subset RB^n$ for $K\in\mathcal{K}^n_o$ and $r,R>0$, then
$$
c_\varphi r B^n\subset \Pi^*_\varphi K\subset 2c_\varphi RB^n.
$$

\begin{lemma}
\label{O*projdefstab}
There exist $\gamma_0\in(0,1]$ depending on $n$ and $\varphi\in\mathcal{C}$
such that  if $\eta\in[0,1)$, $x\in \R^n$ and $K$ is an $o$-symmetric convex body, then
$$
\int_{S^{n-1}}\varphi\left(\frac{x\cdot w}{h_K(w)}\right)\,d V_K(w)\leq
1-\eta\mbox{ \ yields \ }
x\in(1-\gamma_0\cdot \eta)\,\Pi^*_\varphi K.
$$
\end{lemma}
\proof It follows from the linear covariance (\ref{linearmap}) of the polar Orlicz projection body and from John's theorem (see F. John \cite{Joh37}) that we may assume
$$
B^n\subset K\subset \sqrt{n} \, B^n.
$$
Thus  the form of Lemma~2.2 in \cite{LYZ10} above yields
$\Pi^*_\varphi K\subset 2c_\varphi\sqrt{n} \, B^n$.
 
According to (\ref{O*projdef}), there exist   $y\in \partial \Pi^*_\varphi K$ and $\varepsilon\in(0,1)$ such that $x=(1-\varepsilon)y$, and hence if $w\in S^{n-1}$, then
$$
\frac{|y\cdot w|}{h_K(w)}\leq 2c_\varphi\sqrt{n}.
$$
Setting $\gamma_1=\max\{\varphi'(2c_\varphi\sqrt{n}),
-\varphi'(-2c_\varphi\sqrt{n})\}$, we deduce  from
the convexity of $\varphi$ and (\ref{phimonotone}) that if
$t\in[-2c_\varphi\sqrt{n}, 2c_\varphi\sqrt{n}]$, then 
$$
\varphi((1-\varepsilon)t)\geq \varphi(t)-\gamma_1\varepsilon\cdot |t|\geq
\varphi(t)-2c_\varphi\sqrt{n}\cdot\gamma_1\varepsilon.
$$
For $\gamma_2=2c_\varphi\sqrt{n}\cdot\gamma_1$, it follows from (\ref{O*projdef}) that
$$
\int_{S^{n-1}}\varphi\left(\frac{(1-\varepsilon)y\cdot w}{h_K(w)}\right)\,d V_K(w)\geq
\int_{S^{n-1}}\varphi\left(\frac{y\cdot w}{h_K(w)}\right)
-\gamma_2\varepsilon\,d V_K(w)=
1-\gamma_2\varepsilon.
$$
Therefore we may choose $\gamma_0=\min\{1,1/\gamma_2\}$.
\theend

An essential tool to prove Theorem~\ref{OrliczPettyequa}
was the reduction to convex bodies with axial rotational symmetry such that
the shadow boundaries in the directions parallel and orthogonal to the axis are thin.
The core of the argument for Theorem~\ref{OrliczPettystab} is a stability version
 of this reduction, Theorem~\ref{equatorpole}.
To state Theorem~\ref{equatorpole}, we use the following terminology.
We say that a convex body $K$ in $\R^n$ spins around a $u\in S^{n-1}$,
if $K$ is $o$-symmetric, $u\in \partial K$, the axis of rotation of $K$
is $\R u$, and $K\cap u^{\bot}=B^n\cap u^{\bot}$.

\begin{theo}
\label{equatorpole}
Let $K$ be a convex body  in $\R^n$, $n\geq 3$,
such that $\delta_{\rm BM}(K,B^n)\geq\delta\in(0,\delta_0)$,
where $\delta_0>0$ depends on $n$.
Then there exist  $\varepsilon\in(\delta^{24},\delta]$ and a convex body $K'$
spinning around a $u\in S^{n-1}$, such that $K'$ is obtained from
$K$ by a combination of Steiner symmetrizations, linear transformations
and taking limits, and satisfies $\delta_{\rm BM}(K',B^n)\leq\varepsilon$, and
\begin{description}
\item{(i)}
for any $o$-symmetric ellipsoid
$E$ with axial rotational symmetry around $\R u$, one finds
a  ball $x+\varepsilon^2\,B^n\subset {\rm int}(E\Delta K')$
where $|x\cdot u |\leq 1-\varepsilon^2$;
\item{(ii)} $(1-\varepsilon^{32})u+\varepsilon^3 v\not\in K'$ for $v\in S^{n-1}\cap u^\bot$;
\item{(iii)} $\varepsilon^3\,u+(1-\varepsilon^7) v\not\in K'$ for
$v\in S^{n-1}\cap u^\bot$.
\end{description}
\end{theo}

The proof of Theorem~\ref{equatorpole}, being rather technical, is deferred to Section~\ref{secequatorpole}.

As  $\delta_{\rm BM}(K,B^n)\leq \ln n$, Theorem~\ref{OrliczPettystab} follows from
 the following statement. For $\varphi\in\mathcal{C}$, if $K\in\mathcal{K}_o^n$  with
$\delta_{BM}(K,B^n)\geq\delta\in(0,\delta_*)$, then
\begin{equation}
\label{OrliczPettystab0}
\frac{V(\Pi_\varphi^* K)}{V(K)}\leq (1-\gamma\cdot \delta^{792n}\cdot \tilde{\varphi}(\delta^{840}))\,
\frac{V(\Pi_\varphi^* B^n)}{V(B^n)}
\end{equation}
where $\delta_*,\gamma>0$ depend on $n$ and $\varphi$.
In the following the implied constants in $O(\cdot)$ depend on $n$ and $\varphi$.

We always assume that $\delta_*$ in (\ref{OrliczPettystab0}),
and hence $\delta$ and $\varepsilon$, as well, are small
enough to make the argument work. In particular,
$\delta_*\leq\delta_0$ where $\delta_0>0$ is the constant depending $n$
and $\varphi$ of Theorem~\ref{equatorpole}.
It follows from the continuity of the polar Orlicz projection body
that we may also assume the following. If $M$ is a convex body
spinning around a $u\in S^{n-1}$, and $\delta_{\rm BM}(M,B^n)< \delta_*$,
then
\begin{equation}
\label{11}
0.9\, B^n\subset M\subset 1.1\, B^n
\mbox{ \ and \ }0.9\Pi^*_\varphi B^n\subset\Pi^*_\varphi M\subset 1.1
\Pi^*_\varphi\, B^n.
\end{equation}

Let $u_*$ and $\tilde{u}_*$ be orthogonal
unit vectors in $\R^n$, and
let  $K\in\mathcal{K}_o^n$ with $\delta_{\rm BM}(K,B^n)\geq\delta\in(0,\delta_*)$.
According to Theorem~\ref{equatorpole}, there exist $\varepsilon\in(\delta^{24},\delta]$ and a convex body $K'$
spinning around $u_*$ with  $\delta_{\rm BM}(K',B^n)\leq\varepsilon$ and obtained from
$K$ by a combination of Steiner symmetrizations, linear transformations
and taking limits such that
\begin{description}
\item{(i)}
for any $o$-symmetric ellipsoid
$E$ with axial rotational symmetry around $\R u_*$, one finds
a  ball $x+\varepsilon^2\,B^n\subset {\rm int}(E\Delta K')$
where $| x\cdot u_* |\leq 1-\varepsilon^2$;
\item{(ii)} $(1-\varepsilon^{32})u_*+\varepsilon^3 \tilde{u}_*\not\in K'$;
\item{(iii)} $\varepsilon^3\,u_*+(1-\varepsilon^7) \tilde{u}_*\not\in K'$.
\end{description}
It follows from (\ref{linearmap}) and 
 Lemma~\ref{lem14} that
$V(\Pi_\varphi^*K')/V(K')\geq V(\Pi_\varphi^*K)/V(K)$.

We deduce that if $\widetilde{K}$ is a smooth and strictly convex  body
spinning around $u_*$ sufficiently close to $K'$, then
\begin{description}
\item{(a)}
for any $o$-symmetric ellipsoid
$E$ with axial rotational symmetry around $\R u_*$, one finds
a  ball $x+\varepsilon^2\,B^n\subset {\rm int}(E\Delta \widetilde{K})$
where $| x\cdot u_* |\leq 1-\varepsilon^2$;
\item{(b)} $(1-\varepsilon^{32})u_*+\varepsilon^3 \tilde{u}_*\not\in \widetilde{K}$;
\item{(c)} $\varepsilon^3\,u_*+(1-\varepsilon^7) \tilde{u}_*\not\in \widetilde{K}$;
\item{(d)}
$\frac{V(\Pi_\varphi^* \widetilde{K})}{V(\widetilde{K})}\geq
(1-\varepsilon^{33n}\tilde{\varphi}(\varepsilon^{35}))\cdot
\frac{V(\Pi_\varphi^* K)}{V(K)}$;
\item{(e)} $\delta_{\rm BM}(\widetilde{K},B^n)< \delta_*$.
\end{description}

 We  define $v\in S^{n-1}$ by
$$
\lambda_*\,v=\varrho_{\Pi_\varphi^*\widetilde{K}}(u_*)\cdot u_*-
\varrho_{\Pi_\varphi^*\widetilde{K}}(\tilde{u}_*)\cdot \tilde{u}_*
$$
for some $\lambda_*>0$.
It follows from (e) and (\ref{11}) that
\begin{equation}
\label{v*size}
\frac12<\frac{0.9}{\sqrt{0.9^2+1.1^2}}\leq  v\cdot  u_*
\leq \frac{1.1}{\sqrt{0.9^2+1.1^2}}<\frac{\sqrt{3}}2.
\end{equation}
We plan to apply Steiner symmetrization to $\widetilde{K}$
with respect to $v^\bot$, and show that the volume of the polar Orlicz projection
body increases substantially. We consider $v^\bot$ as $\R^{n-1}$,
and set
$$
v^\bot\cap B^n=B^{n-1}.
$$
For $X\subset v^\bot$, the interior of $X$ with respect to the
subspace topology of $v^\bot$ is denoted by ${\rm relint}\,X$.

\begin{figure}
\begin{center}
\includegraphics[width=10em]{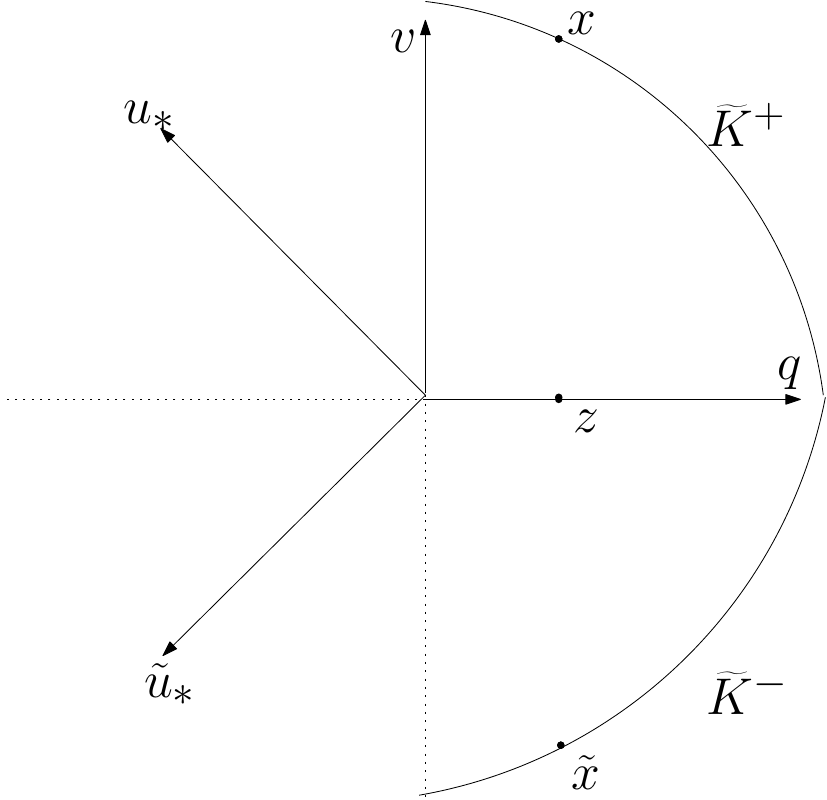}
\end{center}
\caption{ }
\end{figure}

Let $q$ be the unit vector in the line
${\rm lin}\{u_*,\tilde{u}_*\}\cap v^\bot$ satisfying $ q\cdot u_*<0$ (see Figure~1).
We observe that $\Xi_{u_*,\widetilde{K}}=u_*^\bot\cap \partial B^n$
and $\Xi_{\tilde{u}_*,\widetilde{K}}=\tilde{u}_*^\bot\cap \partial \widetilde{K}$, moreover
\begin{eqnarray*}
E_*&=&\pi_{v}(u_*^\bot\cap B^n),\\
\widetilde{K}_*&=&\pi_{v}(\tilde{u}_*^\bot\cap \widetilde{K})
\end{eqnarray*}
are $o$-symmetric, and
have $\R q$ as their axis of rotation inside $v^\bot$.
We define $\theta>0$ by $\theta q\in \partial \widetilde{K}_*$,
and  the linear transform $\Phi:v^\bot\to v^\bot$
by $\Phi(\theta q)=q$,  and $\Phi(y)=y$ for $y\in v^\bot\cap u_*^\bot$.
Thus $\theta\in(\frac12,\frac{\sqrt{3}}2)$ by (\ref{v*size}).
Since $\Phi \widetilde{K}_*$ is congruent to $\tilde{u}_*^\bot\cap \widetilde{K}$,
(a) yields an $(n-1)$-ball
$z'+\varepsilon^2 B^{n-1}\subset {\rm relint}\,\Phi(E_*)\Delta \Phi(\widetilde{K}_*)$
where $0\leq z'\cdot q \leq 1-\varepsilon^2$.
We define $z_*=\Phi^{-1}z'$, and hence
\begin{equation}
\label{x*between}
z_*+\mbox{$\frac{\varepsilon^2}2$}\, B^{n-1}\subset {\rm relint}\,E_*\Delta \widetilde{K}_*.
\end{equation}
Since $ v^\bot\cap u_*^\bot\cap E_*= v^\bot\cap u_*^\bot\cap \widetilde{K}_*$,
we also deduce that
\begin{equation}
\label{x*w*}
\varepsilon^2/2< z_*\cdot q < \theta-(\varepsilon^2/2).
\end{equation}
We write $w_x$ to denote the exterior unit normal
at an $x\in\partial \widetilde{K}$, and  define
\begin{eqnarray*}
\widetilde{K}^+&=&\{x\in\partial\widetilde{K}:\, v\cdot w_x>0
\mbox{ and } q\cdot x>0\};\\
\widetilde{K}^-&=&\{x\in\partial\widetilde{K}:\, v\cdot w_x<0
\mbox{ and } q\cdot x>0\}.
\end{eqnarray*}
It follows from (\ref{x*w*}) that
\begin{equation}
\label{Kproj}
 z_*+\mbox{$\frac{\varepsilon^2}2$}\,B^{n-1}\subset \pi_{v}\widetilde{K}^{\pm}.
\end{equation}
If $z=\pi_{v}x=\pi_{v}\tilde{x}\in  z_*+\frac{\varepsilon^2}4\,B^{n-1}$
for suitable $x\in\widetilde{K}^+$ and $\tilde{x}\in\widetilde{K}^-$,
then $ z+\frac{\varepsilon^2}4\,B^{n-1}\subset \pi_{v}\widetilde{K}^{\pm}$ by (\ref{Kproj}).
We deduce from $\widetilde{K}\subset 1.1B^n$ (compare (\ref{11})) that
$w_x\cdot v, |w_{\tilde{x}}\cdot v|>\varepsilon^2/8$,
and hence (\ref{lrangle}) and (\ref{11}) yield
\begin{equation}
\label{claimfgangle}
0.9\leq\langle f\rangle(z),\langle g\rangle(z)<9\varepsilon^{-2}
\mbox{ \ for $z\in  z_*+\frac{\varepsilon^2}4\,B^{n-1}$}.
\end{equation}

\begin{lemma}
\label{stablemma1}
If $\pi_{v}x=\pi_{v}\tilde{x}\in  z_*+\frac{\varepsilon^2}4\,B^{n-1}$
for $x\in\widetilde{K}^+$ and $\tilde{x}\in\widetilde{K}^-$, then
$$
| u_*\cdot w_x|,\;| \tilde{u}_*\cdot w_{\tilde{x}}|>\varepsilon^{32}/2
\mbox{ \ and \  } (u_*\cdot w_x)\cdot( \tilde{u}_*\cdot w_{\tilde{x}})<0.
$$
\end{lemma}
\proof Since $\widetilde{K}$ is a smooth and strictly convex  body,
and has $\R u_*$ as its axis of rotation,
we have
\begin{eqnarray*}
\Xi^+_{u_*,\widetilde{K}}&=&\{x\in\partial K:\,x\cdot u_*>0\};\\
\Xi^+_{\tilde{u}_*,\widetilde{K}}&=&\{x\in\partial K:\,x\cdot \tilde{u}_*>0\}.
\end{eqnarray*}
It follows from $x\in \widetilde{K}^+$ and $\Xi_{u_*,\widetilde{K}}=u_*^\bot\cap S^{n-1}$ that
\begin{equation}
\label{u*+}
 u_*\cdot w_x>0\mbox{ if and only if }
\pi_{v}x\in{\rm relint}\,\pi_{v}(u_*^\bot\cap B^n)={\rm relint}\,E_*,
\end{equation}
and from $\tilde{x}\in \widetilde{K}^-$ that
\begin{equation}
\label{tildeu*+}
 \tilde{u}_*\cdot w_{\tilde{x}}>0\mbox{ if and only if }
\pi_{v}\tilde{x}\in{\rm relint}\,\pi_{v}(\tilde{u}_*^\bot\cap \widetilde{K})
={\rm relint}\,\widetilde{K}_*.
\end{equation}

We deduce from (\ref{x*between}),  (\ref{u*+}) and (\ref{tildeu*+}) that
\begin{equation}
\label{claim1}
( u_*\cdot w_x)\cdot( \tilde{u}_*\cdot w_{\tilde{x}})<0.
\end{equation}

To have a lower estimate on $| u_*\cdot u_x|$,
we observe that combining (\ref{x*between}) with $\varepsilon^2/4>2\varepsilon^3$
and the fact that $\pi_v$ does not inrease distance yields
$$
x\not \in (u_*^\bot\cap  \partial\widetilde{K})+2\varepsilon^3B^n.
$$
Thus, we conclude from (c) that $\|\pi_{u_*}x\|\leq 1-\varepsilon^7$.
It follows that $(\pi_{u_*}x)+\frac{\varepsilon^7}2\,B^n\subset\widetilde{K}$,
and hence $w_x$ is an exterior normal also to the convex hull at $x$ of this ball and $x$.
As $| u_*\cdot x|\leq 1$, we deduce that
\begin{equation}
\label{claim2}
| u_*\cdot w_x|\geq \varepsilon^7/2.
\end{equation}
Finally we consider $| \tilde{u}_*\cdot w_{\tilde{x}}|$.
Using again (\ref{x*between}), we have
\begin{equation}
\label{tildexe3}
\tilde{x}\not \in (\tilde{u}_*^\bot\cap  \partial\widetilde{K})+2\varepsilon^3B^n.
\end{equation}
In particular $\|\tilde{x}-u_*\|>2\varepsilon^3$ and $\|\tilde{x}-(-u_*)\|>2\varepsilon^3$, and hence (b) implies that
$|\tilde{x}\cdot u_*|< 1-\varepsilon^{32}$. As $\widetilde{K}$ spins around $u_*$, we deduce from (\ref{tildexe3}) that
$(\pi_{\tilde{u}_*}\tilde{x})+\frac{\varepsilon^{32}}2\,B^n\subset\widetilde{K}$.
Thus $w_{\tilde{x}}$ is an exterior normal also to the convex hull at $\tilde{x}$ of this ball and $\tilde{x}$, 
and  hence $| \tilde{u}_*\cdot \tilde{x}|\leq 1$ yields that
\begin{equation}
\label{claim3}
| \tilde{u}_*\cdot w_{\tilde{x}}|\geq \varepsilon^{32}/2.
\end{equation}
Therefore Lemma~\ref{stablemma1} is a consequence
of (\ref{claim1}), (\ref{claim2}) and (\ref{claim3}).
\theend

We continue with the proof
of  Theorem~\ref{OrliczPettystab}.
We use the notation of Lemma~\ref{lem12}. In particular we
write $(z,t)$ to denote $z+tv$ for $z\in\R^{n-1}=v^\bot$ and $t\in\R$,
and  $f$ and $g$ to denote the concave functions
on $\pi_v\widetilde{K}$ such that for
 $z\in{\rm relint}\pi_v\widetilde{K}$, we have $f(z)>-g(z)$,
and $(z,f(z)),(z,-g(z))\in\partial\widetilde{K}$.

We write $\gamma_1,\gamma_2,\ldots$ to denote positive constants depending
 on $n$ and $\varphi$,
and we define
\begin{eqnarray*}
y_*&=&-\pi_v\left(\varrho_{\Pi_\varphi^*\widetilde{K}}(u_*)\cdot u_*\right)=
-\pi_v\left(\varrho_{\Pi_\varphi^*\widetilde{K}}(\tilde{u}_*)\cdot \tilde{u}_*\right),\\
\Psi&=&\left\{\alpha\in S^{n-1}:\, \alpha\cdot v>0\mbox{ \ and  }
\pi_v\left(\varrho_{S_v\Pi*\widetilde{K}}(\alpha)\cdot\alpha\right)
\in y_*+\varepsilon^{33}B^{n-1}\right\}.
\end{eqnarray*}
As $0.9 \Pi_\varphi^*B^n\subset S_v\Pi_\varphi^*\widetilde{K}\subset
1.1\Pi_\varphi^*B^n$
by (e) and (\ref{11}), we have
\begin{equation}
\label{Psisize}
\mathcal{H}(\Psi)>\gamma_1\varepsilon^{33(n-1)}.
\end{equation}
Let
$$
y\in y_*+\varepsilon^{33}B^{n-1},
$$
and  let $(y,t),(y,-s)\in\partial \Pi_\varphi^*\widetilde{K}$
where $-s<t$, and hence $(y,\frac{t+s}2)\in \partial S_v\Pi_\varphi^*\widetilde{K}$.
We define
\begin{eqnarray*}
u&=&\frac{(y,t)}{\|(y,t)\|} \\
\tilde{u}&=&\frac{(y,-s)}{\|(y,-s)\|} \\
\alpha&=&\frac{(y,\frac{t+s}2)}{\|(y,\frac{t+s}2)\|}\in \Psi.
\end{eqnarray*}
It follows from $0.9 \Pi_\varphi^*B^n\subset
 \Pi_\varphi^*\widetilde{K}\subset 1.1 \Pi_\varphi^*B^n$ that
$$
\|u-u_*\|,\|\tilde{u}-\tilde{u}_*\|<\gamma_2\varepsilon^{33}.
$$
Choose $\delta_*$ small enough such that
$\gamma_2\varepsilon^{33}<\varepsilon^{32}/4$.
We deduce from Lemma~\ref{stablemma1} that if
$$
z=\pi_{v}x=\pi_{v}\tilde{x}\in  z_*+(\varepsilon^2/4)\,B^{n-1}
$$
for $x\in\widetilde{K}^+$ and $\tilde{x}\in\widetilde{K}^-$, then
$$
| u\cdot w_x|,\;| \tilde{u}\cdot w_{\tilde{x}}|>\varepsilon^{32}/4
\mbox{ \ and \  }( u\cdot w_x)\cdot( \tilde{u}\cdot w_{\tilde{x}})<0.
$$
Using now $0.9\Pi_\varphi^*B^n\subset \Pi_\varphi^*\widetilde{K}\subset
1.1\Pi_\varphi^*B^n$ and (\ref{lrangle}), we deduce
\begin{equation}
\label{stabcond}
|t-y\cdot \nabla f(z)|,\,|s-y\cdot \nabla g(z)|>\gamma_3\varepsilon^{32}
\mbox{ \ and \  }
(t-y\cdot \nabla f(z))\cdot(s-y\cdot \nabla g(z))<0.
\end{equation}
It follows from (\ref{claimfgangle}) and  (\ref{stabcond}) that we may apply
Lemma~\ref{phiaround0} with
$$
\mbox{$a=t-y\cdot \nabla f(z)$, $b=s-y\cdot \nabla g(z)$,
$\alpha=\langle f\rangle(z)$ and $\beta=\langle g\rangle(z)$.}
$$
 By (\ref{claimfgangle}), (\ref{stabcond}), and since
$\gamma_3\varepsilon^{34}/9>\varepsilon^{35}$,
we may choose $\omega=\varepsilon^{35}$ in Lemma~\ref{phiaround0}, and hence
 (\ref{stabcond}) yields that
\begin{align*}
&\frac12\int_{\pi_v \widetilde{K}}
\varphi\left(\frac{t- y\cdot\nabla f(z)}{\langle f\rangle(z)}\right)\langle f\rangle(z)\,dz
+ \frac12\int_{\pi_v \widetilde{K}}
\varphi\left(\frac{s- y\cdot\nabla g(z)}{\langle g\rangle(z)}\right)\langle g\rangle(z)\,dz\\
&-\int_{\pi_v \widetilde{K}}
\varphi\left(\frac{\frac{t}2+\frac{s}2- \frac{y\cdot\nabla f(z)}2-
\frac{y\cdot\nabla g(z)}2}{\frac{\langle f\rangle(z)}2+\frac{\langle g\rangle(z)}2}\right)
\left(\frac{\langle f\rangle(z)}2+\frac{\langle g\rangle(z)}2\right)\,dz\geq
\gamma_4\,\varepsilon^{-3}\tilde{\varphi}(\varepsilon^{35}).
\end{align*}
Therefore   (\ref{long1}), (\ref{phiconvex}) and (\ref{11}) lead to
\begin{equation}
\label{longstab}
\int_{S^{n-1}}\varphi\left(\frac{(y,\frac12(t+s))\cdot w}{h_{S_v\widetilde{K}}(w)}\right)
\,d V_{S_v\widetilde{K}}(w)\leq
1-\frac{\gamma_4\,\varepsilon^{-3}\tilde{\varphi}(\varepsilon^{35})}{nV(1.1B^n)}.
\end{equation}
We conclude first applying
Lemma~\ref{O*projdefstab}, then the consequence
$0.9\Pi_\varphi^*B^n\subset S_v\Pi_\varphi^*\widetilde{K}
\subset 1.1\Pi_\varphi^* B^n$ of (\ref{11}) that if $w\in \Psi$, then
$$
\varrho_{S_v\Pi_\varphi^*\widetilde{K}}(w)^n\leq
(1-\gamma_5\,\varepsilon^{-3}\tilde{\varphi}(\varepsilon^{35}))^n\cdot
 \varrho_{\Pi_\varphi^*S_v\widetilde{K}}(w)^n
\leq \varrho_{\Pi_\varphi^*S_v\widetilde{K}}(w)^n-\gamma_6\,
\varepsilon^{-3}\tilde{\varphi}(\varepsilon^{35}).
$$
Since $\varrho_{S_v\Pi_\varphi^*\widetilde{K}}(w)\leq
\varrho_{\Pi_\varphi^*S_v\widetilde{K}}(w)$ for any $w\in S^{n-1}$
by Lemma~\ref{lem14}, combining (\ref{radialvol}) and
(\ref{Psisize}) leads to
\begin{eqnarray*}
V(\Pi_\varphi^*\widetilde{K})&=&
V(\Pi_\varphi^*S_v\widetilde{K})\leq V(\Pi_\varphi^*S_v\widetilde{K})-
\gamma_7\,\varepsilon^{33(n-1)-3}\tilde{\varphi}(\varepsilon^{35})\\
&\leq&
(1-\gamma_8\,\varepsilon^{33n-36}\tilde{\varphi}(\varepsilon^{35}))\cdot V(\Pi_\varphi^*S_v\widetilde{K}).
\end{eqnarray*}
We conclude from (d) and Theorem~\ref{PettyOrlicz} that
\begin{eqnarray*}
\frac{V(\Pi_\varphi^* K)}{V(K)}&\leq&
(1-\varepsilon^{33n}\tilde{\varphi}(\varepsilon^{35})))^{-1}
(1-\gamma_8\,\varepsilon^{33n-36}\tilde{\varphi}(\varepsilon^{35})))\cdot
\frac{V(\Pi_\varphi^* \widetilde{K})}{V(\widetilde{K})}\\
&\leq&
(1-\gamma_9\,\varepsilon^{33n}\tilde{\varphi}(\varepsilon^{35})))\cdot
\frac{V(\Pi_\varphi^* B^n)}{V(B^n)},
\end{eqnarray*}
which, in turn, yields  (\ref{OrliczPettystab0}) by $\varepsilon\geq \delta^{24}$.\theend

\section{Proof of Theorem~\ref{OrliczPettystrictstab}}
\label{secOrliczPettystrictstab}

Naturally, we again need a suitable stability version of (\ref{phiconvex}).

\begin{lemma}
\label{phip}
 Let $\varphi\in\mathcal{C}$ be even such that $\varphi''(t)$ is continuous
and positive for $t>0$.
If $a,b,\alpha,\beta,\omega>0$ satisfy
$\omega\leq\frac{a}{\alpha},\frac{b}{\beta}\leq\omega^{-1}$, then
\begin{align*}
&\mbox{$\alpha\varphi\left(\frac{a}{\alpha}\right)+\beta\varphi\left(\frac{b}{\beta}\right)
-(\alpha+\beta)\varphi\left(\frac{a+b}{\alpha+\beta}\right)$}\\
&\geq
\frac{\min\{\varphi''(t):\,t\in(\omega,\omega^{-1})\}\cdot\min\{\alpha^2,\beta^2\}}{2(\alpha+\beta)}\cdot\left(\frac{a}{\alpha}-\frac{b}{\beta}\right)^2.
\end{align*}
\end{lemma}
\proof The Taylor formula around $\frac{a+b}{\alpha+\beta}$
yields the estimate. \theend

Given Theorem~\ref{OrliczPettystab}, what we need to consider are translates of
a convex body that are close to the unit ball.

\begin{lemma}
\label{Ktranslate}
Let $\varphi\in\mathcal{C}$ be even such that $\varphi''(t)$ is continuous
and positive for $t>0$. There exist $\varepsilon_0,\gamma>0$
depending on $n$ and $\varphi$ such that if
 $\|\theta\|\geq\varepsilon^{1/3}$ and
$B^n\subset K-\theta\subset(1+\varepsilon)B^n$
for $K\in\mathcal{K}_0^n$, $\varepsilon\in (0,\varepsilon_0)$
and $\theta\in\R^n$, then
$$
\frac{V(\Pi_\varphi^* K)}{V(K)}<(1-\gamma \varepsilon^{\frac23})\,
\frac{V(\Pi_\varphi^*B^n)}{V(B^n)}.
$$
\end{lemma}
\proof We write $\sigma$ to denote the reflection through $\theta^\bot$.
Possibly after applying Schwarz rounding with respect to $v=\theta/\|\theta\|$ (compare
Lemma~\ref{lem14}), we may assume that $\R v$ is the axis of rotation of $K$.
It follows that $\Pi_\varphi^* K$ also has $\R v$ as its axis of rotation. Since $\varphi$ is  even, we deduce that $\Pi_\varphi^* K$ is $o$-symmetric, therefore $\Pi_\varphi^* K$ is symmetric with respect to $\sigma$. We may also assume that $K$ is smooth, and we write
$w_x$ to denote the unique exterior unit normal at $x\in\partial K$.

We write $\gamma_1,\gamma_2,\ldots$ to denote positive constants
depending on $n$ and $\varphi$. In addition the implied constant
in $O(\cdot)$ depends also only on $n$ and $\varphi$.
As $K\subset 3B^n$,
Lemma~2.2 by  E. Lutwak, D. Yang, G. Zhang \cite{LYZ10} yields
\begin{equation}
\label{polarinball}
\Pi_\varphi^* K\subset \gamma_1B^n.
\end{equation}
Since $\Pi_\varphi^* K$ is $o$-symmetric and $\Pi_\varphi^* K\subset \gamma_1B^n$,
there exists $\gamma_2>0$ depending on $n$ and $\varphi$, such that
 if $h_{\Pi_\varphi^* K}(u)\leq \gamma_2$
for some $u\in S^{n-1}$, then  $V(\Pi_\varphi^* K)< \frac12\,V(\Pi_\varphi^* B^n)$. In particular,
Lemma~\ref{Ktranslate} readily holds in this case. Therefore
we may assume that
\begin{equation}
\label{ballinpolar}
\gamma_2\,B^n\subset \Pi_\varphi^* K.
\end{equation}

We set $\R^{n-1}=v^\bot$ and $B^{n-1}=v^\bot\cap B^n$,
and write the point $y+tv$ of $\R^n$ with $y\in\R^{n-1}$
and $t\in \R$ in the form $(y,t)$.
In addition, let $f,g$ be the concave functions on $\pi_vK$ satisfying
$$
K=\{(y,t):\,y\in\pi_vK\mbox{ and }-g(y)\leq t\leq f(y)\}.
$$
We consider
\begin{align*}
\Xi&=\mbox{$\frac35\,B^{n-1}\backslash \frac12\,B^{n-1}$}\\
\Psi&=\mbox{$\{(y,t)/\|(y,t)\|\in S^{n-1}:\,y\in\frac{3\gamma_2}5\,B^{n-1},
\mbox{ }t>0\mbox{ and }(y,t)\in\partial \Pi_\varphi^* K\}.$}
\end{align*}
It follows that
\begin{equation}
\label{Psi2}
\mathcal{H}(\Psi)\geq \gamma_3.
\end{equation}

For $y\in\frac{3\gamma_2}5\,B^{n-1}$ and $z\in\Xi$,
 let $t>0$ such that $(y,t)\in\partial \Pi_\varphi^* K$, and hence
$(y,-t)\in\partial \Pi_\varphi^* K$ since  $\sigma(\Pi_\varphi^* K)=\Pi_\varphi^* K$.
We plan to apply Lemma~\ref{phip} with
\begin{equation}
\label{abdef2}
\mbox{$a=t- y\cdot\nabla f(z)$, $b=t- y\cdot\nabla g(z)$, $\alpha=\langle f\rangle(z)$
and $\beta=\langle g\rangle(z)$. }
\end{equation}
Let $x,\tilde{x}\in\partial K$,
and let $x',\tilde{x}'\in\partial (\theta+B^n)$ be defined in a way such that
$\pi_vx=\pi_v\tilde{x}=\pi_vx'=\pi_v'\tilde{x}'=z$, $(x-\tilde{x})\cdot v>0$
and $(x'-\tilde{x}')\cdot v>0$. We  observe that
$\sigma(\tilde{x}'-\theta)=x'-\theta$.
The condition $z\in\Xi$ yields that
\begin{equation}
\label{vx'}
\frac45\leq v\cdot (x'-\theta)=-v\cdot (\tilde{x}'-\theta)\leq \frac{\sqrt{3}}2<0.9.
\end{equation}
Since the angles between $v$ and both $(y,t)$
and $x'-\theta$ are at most $\gamma_4={\rm arcsin}\,\frac35$,
 and $\cos2\gamma_4=\frac{7}{25}$,  we deduce from (\ref{polarinball}) and (\ref{ballinpolar}) that
\begin{equation}
\label{ytx'}
\frac{7\gamma_2}{25}\leq (y,t)\cdot (x'-\theta)=(y,-t)\cdot (\tilde{x}'-\theta)
\leq \gamma_1.
\end{equation}
To compare $x'-\theta$ and $w_x$, we observe that the tangent planes
to $\theta+B^n$ at both $x'$ and $\theta+w_x$ separate $x$ and $\theta+B^n$. Since
$\|x-\theta\|\leq 1+\varepsilon$, such points on
$\theta+S^{n-1}$ are contained in a cap cut off by a hyperplane of distance at least
 $(1+\varepsilon)^{-1}$ from $\theta$, and the diameter of the cap
is at most $2\sqrt{1-(1+\varepsilon)^{-2}}<4\varepsilon^{\frac12}$. Therefore
\begin{equation}
\label{wx'}
\|w_x-(x'-\theta)\|<4\varepsilon^{\frac12}
\mbox{ \ and \ }\|w_{\tilde{x}}-(\tilde{x}'-\theta)\|<4\varepsilon^{\frac12}.
\end{equation}
From (\ref{lrangle}), (\ref{abdef2}), (\ref{ytx'}) and (\ref{wx'}), we deduce that
\begin{eqnarray}
\label{aper}
\frac{a}{\alpha}&=&\left(1+O(\varepsilon^{\frac12})\right)\cdot
\frac{(y,t)\cdot (x'-\theta)}{h_K(w_x)};\\
\label{bper}
\frac{b}{\beta}&=&\left(1+O(\varepsilon^{\frac12})\right)\cdot
\frac{(y,t)\cdot (x'-\theta)}{h_K(w_{\tilde{x}})}.
\end{eqnarray}
We have $\theta\cdot w+1\leq h_K(w)\leq \theta\cdot w+1+\varepsilon$
for any $w\in S^{n-1}$, and $\|\theta\|<1+\varepsilon$ by $o\in{\rm int}\,K$.
Therefore (\ref{vx'}),  (\ref{wx'}) and the
condition $\|\theta\|\geq\varepsilon^{1/3}$ yield
\begin{equation}
\label{hwx}
1+\mbox{$\frac35$}\,\varepsilon^{1/3}<h_K(w_x)<1.9
\mbox{ \ and \ }
0.1<h_K(w_{\tilde{x}})<1-\mbox{$\frac35$}\,\varepsilon^{1/3}
\end{equation}
provided that $\varepsilon_0>0$ is suitably small. We deduce from
(\ref{ytx'}), (\ref{aper}), (\ref{bper}) and (\ref{hwx}) that
there exist $\omega,\gamma_5>0$ depending on $n$ and $\varphi$ such that
\begin{eqnarray}
\label{abest}
\frac{b}{\beta}-\frac{a}{\alpha}&>&\gamma_5\varepsilon^{\frac13};\\
\label{abdiff}
\omega<\frac{a}{\alpha}&<&\frac{b}{\beta}<\omega^{-1}.
\end{eqnarray}
In addition, (\ref{lrangle}), (\ref{vx'}), (\ref{wx'}) and (\ref{hwx}) yield that
\begin{equation}
\label{alphabetaest}
\gamma_6<\alpha,\beta<\gamma_7.
\end{equation}
We conclude from Lemma~\ref{phip} the estimate
\begin{align}
\nonumber
&\frac12\int_{\pi_v K}
\varphi\left(\frac{t- y\cdot\nabla f(z)}{\langle f\rangle(z)}\right)\langle f\rangle(z)\,dz
+ \frac12\int_{\pi_v K}
\varphi\left(\frac{t- y\cdot\nabla g(z)}{\langle g\rangle(z)}\right)\langle g\rangle(z)\,dz\\
\label{zinXi}
&-\int_{\pi_v K}
\varphi\left(\frac{t- \frac{y\cdot\nabla f(z)}2-\frac{y\cdot\nabla g(z)}2}{\frac{\nabla f(z)}2+\frac{\nabla g(z)}2}\right)
\left(\frac{\nabla f(z)}2+\frac{\nabla g(z)}2\right)\,dz
>\gamma_8\varepsilon^{\frac23}.
\end{align}
Since (\ref{zinXi}) holds for any $z\in\Xi$, and
 $S_v\Pi_\varphi^*K=\Pi_\varphi^*K$, we deduce from (\ref{long1}) and  (\ref{phiconvex})
 that
\begin{equation}
\label{long100}
\int_{S^{n-1}}\varphi\left(\frac{(y,t)\cdot w}{h_{S_vK}(w)}\right)\,d V_{S_vK}(w)
<1-\gamma_9\varepsilon^{\frac23}.
\end{equation}
Now we have (\ref{long100}) for all $y\in\frac{3\gamma_2}5\,B^{n-1}$,
and hence
$$
\varrho_{\Pi_\varphi^*K}(u)<
(1-\gamma_{10}\varepsilon^{\frac23})\varrho_{\Pi_\varphi^*S_vK}(u)
$$
for $u\in\Psi\subset S^{n-1}$ according to Lemma~\ref{O*projdefstab},
where $\mathcal{H}(\Psi)\geq \gamma_3$ by (\ref{Psi2}).
Therefore combining Lemma~\ref{lem14}, (\ref{radialvol}), (\ref{polarinball}) and (\ref{ballinpolar})  yields  Lemma~\ref{Ktranslate}.
\theend

Theorem~\ref{OrliczPettystrictstab} follows from the following statement.
For $\varphi\in\mathcal{C}$, there exist $\eta_0,\gamma>0$ depending only on $n$ and $\varphi$  such that if
$K\in\mathcal{K}^n_o$, $\eta\in(0,\eta_0)$, and
\begin{equation}
\label{pcond}
K\not\subset (1+\eta)E \mbox{ \ for any $o$-symmetric
ellipsoid $E\subset K$},
\end{equation}
 then
\begin{equation}
\label{OrliczPettystrictstab0}
\frac{V(\Pi^*_\varphi K)}{V(K)}\leq\left(1-\gamma\cdot \eta^{2376n}\cdot\varphi(\eta^{2520})\right)
\cdot\frac{V(\Pi^*_\varphi B^n)}{V(B^n)}.
\end{equation}

If $\delta_{BM}(K,B^n)>\eta^3/108$, then
Theorem~\ref{OrliczPettystab} yields (\ref{OrliczPettystrictstab0}).
Therefore we assume that
$\delta_{BM}(K,B^n)\leq\eta^3/108$. In particular,
we may assume that
$\theta+B^n\subset K$ for some $\theta\in\R^n$, and $K$ is contained in a ball of radius
$1+\frac{\eta^3}{54}$. It follows that
$$
\theta+B^n\subset K\subset \theta+\left(1+\mbox{$\frac{\eta^3}{27}$}\right)\,B^n.
$$
We deduce from (\ref{pcond}) that
$\frac{1+\|\theta\|+\frac{\eta^3}{27}}{1-\|\theta\|}>1+\eta$, and
hence $\|\theta\|>\eta/3$. Therefore we may apply Lemma~\ref{Ktranslate}
with $\varepsilon=\frac{\eta^3}{27}$, which, in turn, completes the
proof of (\ref{OrliczPettystrictstab0}).
\theend

\section{Class reduction based on Steiner symmetrization}
\label{secequatorpole}

In this section, we prove Theorem~\ref{equatorpole}.
Let
$$
\mbox{  $u\in S^{n-1}$
\ \ and  \ \  $v\in S^{n-1}\cap u^\bot$.}
$$
Recall that a convex body $K$ in $\R^n$ spins around  $u$,
if $K$ is $o$-symmetric, $u\in \partial K,$ the axis of rotation of $K$
is $\R u$, and $K\cap u^{\bot}=B^n\cap u^{\bot}$.
In this case, we call $\pm u$ the poles of $K$, and $\partial K\cap u^{\bot}\subset S^{n-1}$
the equator of $K$.
We show that  to have a stability
version of the Orlicz-Petty projection inequality,
we may assume that $K$ is an $o$-symmetric convex body with axial rotational
symmetry such that the boundary sufficiently bends near the
equator and the poles.

We prepare the proof of Theorem~\ref{equatorpole}
by a series of Lemmas.
First of all, one may assume
that $K$ is an $o$-symmetric convex body with axial rotational
symmetry because of the following.

\begin{lemma}
\label{rounding}
For any $n\geq 2$ there exists $\gamma>0$ depending
only on $n$, such that
if $K$ is a convex body  in $\R^n$ such that
$\delta_{\rm BM}(K,B^n)\geq \varepsilon\in(0,1)$, then one can find an
$o$-symmetric convex body $C$ with
axial rotational symmetry and
$\delta_{\rm BM}(C,B^n)= \gamma\varepsilon^2$
that is obtained from $K$ using Steiner symmetrizations,
linear transformations and taking limits.
\end{lemma}
{\bf Remark: } If $K$ is $o$-symmetric, then
$\delta_{\rm BM}(C,B^n)=\gamma\varepsilon$ is possible.\\
\proof According to Theorem~1.4 in \cite{Bor}
there is an $o$-symmetric convex body $C$ with
axial rotational symmetry
that is obtained from $K$ using Steiner symmetrizations,
linear transformations and taking limits, and that satisfies
$\delta_{\rm BM}(C_0,B^n)\geq \gamma\varepsilon^2$.
We note that in Theorem~1.4, it is stated that affine transformations
are needed. But translations  are only used to translate $K$ at the beginning
by $-\sigma_K$ where $\sigma_K$ is the centroid of $K$.
If we perform all Steiner symmetrizations
in the proof  of Theorem~1.4  in \cite{Bor} through the same
hyperplanes containing the origin, then even
without the translation at the beginning, the convex body $C_0$
will still be $o$-symmetric.

We may assume that $\delta_{\rm BM}(C_0,B^n)> \gamma\varepsilon^2$,
otherwise we are done.
Since some sequence of Steiner symmetrizations subsequently applied to
$C_0$ converges to a Euclidean ball $B_0$ of volume $V(C_0)$, there is
a sequence $\{C_m\}$, $m=0,1,2,\ldots$ of $o$-symmetric convex
bodies tending to $B_0$ such that $C_m$, $m>0$,
is a Schwarz rounding of $C_{m-1}$
with respect to some $w_m\in S^{n-1}$. In particular, there is
$m\geq 0$ such that
$\delta_{\rm BM}(C_m,B^n)> \gamma\varepsilon^2$ and
$\delta_{\rm BM}(C_{m+1},B^n)\leq \gamma\varepsilon^2$.

For $w\in S^{n-1}$, let $M_w$ be the Schwarz rounding of $C_m$
with respect to $\R w$. Then $\delta_{\rm BM}(M_w,B^n)$
is a continuous function of $w$. Since $C_m=M_{w_m}$
and $C_{m+1}=M_{w_{m+1}}$,
there is a $w\in S^{n-1}$ with $\delta_{\rm BM}(M_w,B^n)=\gamma\varepsilon^2$.

 If $K$ is $o$-symmetric, then Theorem~1.4  in \cite{Bor}
states that $\delta_{\rm BM}(C_0,B^n)\geq \gamma\varepsilon$,
and hence the argument above gives
$\delta_{\rm BM}(C,B^n)=\gamma\varepsilon$.
\theend

In order to obtain a stability version of the Orlicz-Petty projection inequality for an
$o$-symmetric convex body $K$ with axial rotational
symmetry, it is hard to deal with $K$ if it is close to be flat at the poles,
or close to be ruled near the equator. In these
cases, we apply an extra Schwarz rounding. The precise
statements are the subjects of Lemma~\ref{equator} and
Proposition~\ref{polefinal}.
 For $w\in S^{n-1}$ and $t\in\R$, we recall that
$$
H(w,t)=w^{\bot}+tw.
$$
The next observation considers the shape of a convex body with axial rotational symmetry
near the equator.

\begin{lemma}
\label{equator}
There exist $\tau_1,\tau_2>0$ depending
on $n$ with the following properties.
 If  $t\in(0,\frac13)$, the  convex body $K$ in $\R^n$ spins around $u$, and
$$
\tau_1\sqrt{t}\,u+(1-t)v\in K,
$$
then $\delta_{\rm BM}(K',B^n)\geq \tau_2 t$
for the Schwarz rounding $K'$ of $K$ around $\R v$.
\end{lemma}
\proof  Let $E_0$ be the  $o$-symmetric ellipsoid with
axial rotational symmetry around $\R v$
such that $v\in\partial E_0$, and $\mathcal{H}\left(E_0\cap v^\bot\right)=2\kappa_{n-1}$.
For any $s\in(0,\frac23)$, we have
\begin{equation}
\label{E0far}
\gamma_1\sqrt{s}\,u+(1-s)v\not\in (1+\tau_2 s)E_0
\end{equation}
for suitable $\gamma_1>0$ and $\tau_2\in(0,1)$ depending only on $n$.
We define $\tau_1$
by the equation 
$$
(\tau_1\kappa_{n-2}/\kappa_{n-1})^{\frac1{n-1}}=\gamma_1\sqrt{2}.
$$

 Let $E\subset K'$ be
an $o$-symmetric ellipsoid with axial rotational symmetry around $\R v$
such that  $K'\subset \lambda E$, where $\ln\lambda=\delta_{\rm BM}(K',B^n)$.
It follows from the normalization of $K$ that $\mathcal{H}\left(K\cap v^\bot\right)\leq 2\kappa_{n-1}$, thus $E\subset E_0$.

If $\tau_1\sqrt{t}\,u+(1-t)v\in K$ for $t\in(0,\frac13)$, then $\tau_1\sqrt{t}\,u+(1-t)(u^\bot\cap B^n)\subset K$ 
and $\sqrt{t(2-3t)}>\sqrt{t}$ yield that 
$$
\tau_1\sqrt{t}\,u+(1-2t)v+\sqrt{t}(u^\bot\cap v^\bot\cap B^n)\subset K.
$$
Since $H(v,1-2t)\cap K$ contains an $(n-1)$-dimensional cylinder whose height is $\tau_1\sqrt{t}$, and whose base has radius
$\sqrt{t}$, we have
$$
\mathcal{H}\left(H(v,1-2t)\cap K'\right)=\mathcal{H}\left(H(v,1-2t)\cap K\right)\geq
\tau_1\kappa_{n-2}t^{\frac{n-1}2}.
$$
In particular
$$
\gamma_1\sqrt{2t}+(1-2t)v=(\tau_1\kappa_{n-2}/\kappa_{n-1})^{\frac1{n-1}}\sqrt{t}\,u+(1-2t)v\in K'.
$$
We conclude from (\ref{E0far}), that $\lambda> 1+\tau_2 2 t$, and hence
$\delta_{\rm BM}(K',B^n)> \tau_2  t$.
\theend

Now we consider the shape of a convex body with axial rotational symmetry
near the poles. To test whether a convex body is "flat"
near the  poles, we will use the following statement.

\begin{lemma}
\label{polebasic}
There exist $\delta_0,\tau_0,\tau\in(0,1)$ depending
on $n$ with the following property.
 Let $\delta\in(0,\delta_0)$, $t\in(0,\tau_0\delta)$, and let a convex body $K$ with
 $\delta=\delta_{\rm BM}(K,B^n)$ spin around $u$. If an
$o$-symmetric ellipsoid
$E$ with axial rotational symmetry around $\R u$ satisfies that
$E\Delta K$ contains no ball of the form $x+t\,B^n$
with $| x\cdot u |\leq 1-t$, then
\begin{description}
\item{(i)} $K\subset (1+\tau t)E$;
\item{(ii)} assuming $| x \cdot u |\leq 1-4t$,
$x\in\partial E$ implies $(x+3 tB^n)\cap K\neq \emptyset$,
and  $x\in\partial K$ implies $(x+3 tB^n)\cap E\neq \emptyset$;
\item{(iii)} $\theta t\in\partial E$ where $1+\frac12\,\delta\leq \theta\leq 1+\tau\delta$.
\end{description}
\end{lemma}
\proof We write $\gamma_1,\gamma_2,\ldots$ to denote
positive constants depending only on $n$.

For an $x\in\R^n$ with $| x \cdot u |\leq 1-4t$,
we may assume that $x\cdot u\geq 0$.
Let $v\in u^\bot$ such that $x\cdot v\geq 0$ and $x\in{\rm lin}\{u,v\}$.
Since $x+3tB^n$  contains  $x-tu-tv+tB^n$,
we deduce (ii) from the assumptions on $E$ and $K$.

As $K$ spins around $u$, and  $\delta_{\rm BM}(K,B^n)=\delta$, we have
$$
(1/2)B^n\subset (1-\gamma_1\delta) B^n
\subset K\subset (1+\gamma_2\delta)B^n.
$$
This combined with (ii) implies (i). In addition we deduce from (ii) that
$$
(1-\gamma_3t)K \subset \{x\in E:\,| x\cdot u |\leq 1-7t\}
\subset (1+\gamma_4t)K,
$$
which in turn  yields that if $\theta t\in\partial E$ for $\theta>0$, then
$$
\delta=\delta_{\rm BM}(K,B^n)\leq \ln\left[
 (1-\gamma_3t)^{-1}\cdot \theta(1-7t)^{-1} (1+\gamma_4t)\right]
\leq  \ln\theta +\gamma_5t.
$$
Therefore assuming  $t<(2\gamma_5)^{-1}\delta$, we have
$\theta\geq 1+\frac{\delta}2$.
\theend

\begin{coro}
\label{pole0}
There exist $\delta_0,\tau_0\in(0,1)$ depending
on $n$ with the following property.
 Let $\delta\in(0,\delta_0)$, $t\in(0,\tau_0\delta)$, and let a convex body $K$ with
 $\delta=\delta_{\rm BM}(K,B^n)$ spin around $u$. If an
$o$-symmetric ellipsoid
$E$ with axial rotational symmetry around $\R u$ satisfies that
$E\Delta K$ contains no ball of the form $x+t\,B^n$
with $| x\cdot u |\leq 1-t$, then
$$
(1-7t)u+(\sqrt{\delta}/4)v\in K.
$$
\end{coro}
\proof By Lemma~\ref{polebasic} (iii),
we have $\theta u\in\partial E$ where $\theta>1+\frac12\,\delta$.
It follows that
$$
\sqrt{1-\frac{(1-4t)^2}{\theta^2}}>
\sqrt{1-\frac{1}{1+\delta}}>\sqrt{\delta}/2,
$$
and hence
$$
w=(1-4t)u+(\sqrt{\delta}/2) v\in E.
$$
Thus, we obtain Corollary~\ref{pole0} from  Lemma~\ref{polebasic} (ii).
\theend

If a convex body  with axial rotational symmetry is "too flat" around the poles
then we modify it in the following way.

\begin{lemma}
\label{pole}
 If  $\varepsilon\in(0,\varepsilon_0)$ for  $\varepsilon_0\in(0,1)$ depending
on $n$, and $K$ is a convex body  with
 $\delta_{\rm BM}(K,B^n)=\varepsilon$ spinning around $u$, then
there exists a convex body $K'$ that spins around $u$,
and is obtained from $K$ by combining  linear
transformations and one Schwarz rounding, such that
for any $o$-symmetric ellipsoid
$E$ with axial rotational symmetry around $\R u$, one finds
a  ball of the form $x+4\varepsilon^2\,B^n$ in
$E\Delta K'$, where $| x\cdot u |\leq 1-4\varepsilon^2$.
\end{lemma}
\proof In the following the implied constants in $O(\cdot)$ depend only
on $n$, and we write $\gamma_1,\gamma_2,\ldots$ to denote
positive constants depending only on $n$. We assume that
$\varepsilon_0$  depends only on $n$ and is small enough
to make the argument below work.

If for any $o$-symmetric ellipsoid
$E$ with axial rotational symmetry around $\R u$, one finds
a  ball of the form $x+\varepsilon^{3/2}\,B^n$ in $E\Delta K$
where $|( x\cdot u)|\leq 1-\varepsilon^{3/2}$,
then we are done. Therefore let us assume that this is not the case, and hence there exists  an $o$-symmetric ellipsoid $E_0$
 with axial rotational symmetry around $\R u$ satisfying that
$E_0\Delta K$ contains no ball of the form $x+\varepsilon^{3/2}\,B^n$
with $| x\cdot u |\leq 1-\varepsilon^{3/2}$.
Let $u$ be part of an orthonormal
basis for $\R^n$, let $\Phi$ be the diagonal matrix
that maps $E_0$ into $B^n$, and let $K_0=\Phi K$.

By Lemma~\ref{polebasic} (iii) applied to $K$ and $E_0$,
we have $\theta u\in\partial E_0$ where
 $1+\frac12\,\varepsilon<\theta<1+\gamma_1\,\varepsilon$, and hence
$$
(1-s)u\in\partial K_0,\mbox{ \ \  where
$\frac14\,\varepsilon<s<\gamma_2\,\varepsilon $.}
$$
In addition,  Lemma~\ref{polebasic} (i) and (ii) yield
\begin{eqnarray*}
K&\subset&\left (1+\gamma_3\varepsilon^{3/2}\right)E_0,\\
(x+3 \varepsilon^{\frac32}B^n)\cap K&\neq& \emptyset
\mbox{ \ for all $x\in \partial E_0$ with $| x\cdot u |\leq 1-4\varepsilon^{3/2}$}.
\end{eqnarray*}
Thus, we deduce that
\begin{align}
\label{K0Bin}
&K_0\subset\left (1+\gamma_3\varepsilon^{\frac32}\right)B^n,\\
\label{K0Bout}
&(x+4 \varepsilon^{\frac32}B^n)\cap K\neq \emptyset
\mbox{ \ for all $x\in S^{n-1}$ with $| x\cdot u |\leq 1-s-4\varepsilon^{\frac32}$}.
\end{align}
Since $\frac14\,\varepsilon<s<\gamma_2\,\varepsilon $ implies
$$
\sqrt{\left (1+\gamma_3\varepsilon^{\frac32}\right)^2-
\left( 1-s-8\varepsilon^{\frac32}\right)^2}>\sqrt{2s}-\gamma_5\varepsilon,
$$
we deduce from (\ref{K0Bout})  that
\begin{equation}
\label{K0Bout0}
(1-s-8\varepsilon^{\frac32})u+\left(\sqrt{2s}-\gamma_5\varepsilon\right)v\in K_0.
\end{equation}

We plan to apply Schwarz rounding of $K_0$ with respect to
$\R u'$, where
$$
u'=\sqrt{1-s}\,u+\sqrt{s}\,v.
$$
It follows from $\sqrt{1-s}=1-\frac12\,s+O(s^2)$, (\ref{K0Bin})
and (\ref{K0Bout0}) that
\begin{equation}
\label{hest}
1-\left(\mbox{$\frac32$}-\sqrt{2}\right) s-\gamma_6\varepsilon^{\frac32}
\leq h_{K_0}(u')\leq 1+\gamma_3\varepsilon^{\frac32}.
\end{equation}
Next let
$$
\varepsilon^{\frac32}/2<p<2\varepsilon^{\frac32},
$$
let $w$ be of the form $w=(1-s)\,u+t\,v$ with
 $ w\cdot u'=h_{K_0}(u')-p$, and let
$z=(h_{K_0}(u')-p)u'$. In addition, let $\varrho$ be the radius of
$$
G=H(u',h_{K_0}(u')-p)\cap(1+\gamma_3\varepsilon^{\frac32})B^n.
$$
As $H(u',h_{K_0}(u')-p)$ cuts of a cap of depth at most
$(\frac32-\sqrt{2}+O(\varepsilon^{\frac12}))\cdot s$ from
$(1+\gamma_3\varepsilon^{\frac32})B^n$ by (\ref{hest}), and
$\frac32-\sqrt{2}=\frac12(\sqrt{2}-1)^2$, we have
$$
\varrho\leq \left((\sqrt{2}-1)+O(\varepsilon^{\frac12})\right)\sqrt{s}.
$$
In addition, for $y=\sqrt{1-s}\,u'$ (collinear with $w$ and $(1-s)u$),
we have
$$
\|y-z\|\geq \left(\sqrt{2}-1-O(\varepsilon^{\frac12})\right)s,
$$
therefore
$$
\|w-z\|=\frac{\sqrt{1-s}}{\sqrt{s}}\,\|y-z\|
\geq \left(\sqrt{2}-1-O(\varepsilon^{\frac12})\right)\sqrt{s}.
$$
Now $H(u,1-s)$ cuts of a cap of depth
$$
\varrho-\|w-z\|\leq O(\varepsilon^{\frac12})\sqrt{s}=O(\varepsilon)
$$
from $G$, and this cap contains $H(u',h_{K_0}(u')-p)\cap K_0$.
We deduce that
$$
\mathcal{H}\left(H(u',h_{K_0}(u')-p)\cap K_0\right)\leq
O(\varepsilon)(\varepsilon\,\varrho)^{\frac{n-2}2}
\leq O(\varepsilon^{\frac14}) \varepsilon^{\frac{3(n-1)}4}.
$$
Let $K_1$ be the Schwarz rounding of $K_0$ around $\R u'$,
and let $K'$ be the convex body spinning around $u$
that is the image of $K_1$ by a linear transformation that
maps $h_{K_1}(u')u'$ into $u$, and $K_1\cap u'^\bot$ into $B^n\cap u^\bot$.
Thus $K'$ satisfies
$$
\mathcal{H}\left(H(u,1-\varepsilon^{\frac32})\cap K'\right)\leq
 O(\varepsilon^{\frac14}) \varepsilon^{\frac{3(n-1)}4}.
$$
We conclude that $\delta_{\rm BM}(K',B^n)\geq  \gamma_6\varepsilon^{\frac{3}2}$
on the one hand, and
\begin{equation}
\label{closetoaxis}
(1-\varepsilon^{\frac32})u+\gamma_7\varepsilon^{\frac1{4(n-1)}}\cdot \varepsilon^{\frac{3}4}v\not\in K'
\end{equation}
on the other hand.

Next we suppose that there exists
some  $o$-symmetric ellipsoid
$E$ with axial rotational symmetry around $\R u$, such that
no  ball of the form $x+4\varepsilon^2\,B^n$
with  $| x\cdot u |\leq 1-4\varepsilon^2$ is contained in
$E\Delta K'$. By Lemma~\ref{pole0} and
$\delta_{\rm BM}(K',B^n)\geq  \gamma_6\varepsilon^{\frac{3}2}$, we have
\begin{equation}
\label{fartoaxis}
(1-28\varepsilon^2)u+ \gamma_8\varepsilon^{\frac34} \,v\in K'.
\end{equation}
If $\varepsilon_0$ is small enough, then (\ref{closetoaxis})
contradicts  (\ref{fartoaxis}), completing
the proof of Lemma~\ref{pole}.
\theend

Next, strengthening  Lemma~\ref{pole}, we are even more specific
about the shape of the $o$-symmetric convex body
with axial rotational symmetry near the poles.

\begin{prop}
\label{polefinal}
 If  $\varepsilon\in(0,\varepsilon_0)$
for  $\varepsilon_0\in(0,1)$ depending
on $n$, and $K$ is a convex body  spinning around $u$ such that
 $\delta_{\rm BM}(K,B^n)=\varepsilon$, then
there exists a convex body $K'$ that spins around $u$,
and is obtained from $K$ by combining  linear
transformations and two Schwarz roundings, such that
\begin{description}
\item{(i)}
for any $o$-symmetric ellipsoid
$E$ with axial rotational symmetry around $\R u$, one finds
a  ball $x+2\varepsilon^2\,B^n\subset E\Delta K'$
where $| x\cdot u |\leq 1-2\varepsilon^2$;
\item{(ii)} $(1-\varepsilon^{32})u+\varepsilon^3 v\not\in K'$.
\end{description}
\end{prop}

\proof   In the following the implied constants in $O(\cdot)$ depend only on $n$. We assume that
$\varepsilon_0$  depends only on $n$ and is small enough
to make the argument below work.

According to Lemma~\ref{pole}, there exists
a convex body  $K_0$ that spins around $u$,
and is obtained from $K$ by combining  linear
transformations and a Schwarz rounding, such that
for any $o$-symmetric ellipsoid
$E$ with axial rotational symmetry around $\R u$  and
$E\cap u^\bot=B^n\cap u^\bot$, one finds
a  ball of the form $x+2\varepsilon^4\,B^n$ in
$E\Delta K_0$
where $| x\cdot u |\leq 1-4\varepsilon^2$.
If  $(1-\varepsilon^{32})u+\varepsilon^3 v\not\in K_0$,
then we may take $K'=K_0$. Therefore we assume
that
\begin{equation}
\label{e14cond}
(1-\varepsilon^{32})u+\varepsilon^3 v\in K_0.
\end{equation}
To obtain $K'$, first we apply Schwarz rounding around
$\R u'$ for the unit vector
$$
\tilde{u}=\sqrt{1-\varepsilon^{32}}\,u+\varepsilon^{16}\, v
$$
to get  a convex body $\widetilde{K}$. Then we set $K'=\widetilde{\Phi}\widetilde{K}$
where $\widetilde{\Phi}$ is a linear transform that maps
$h_{\widetilde{K}}(\tilde{u})\tilde{u}=h_{K_0}(\tilde{u})\tilde{u}$ into $u$,
and $\widetilde{K}\cap\tilde{u}^\bot$ into $B^n\cap u^\bot$.

Since  $\delta_{\rm BM}(K_0,B^n)\leq\varepsilon$, we have
\begin{equation}
\label{K0K1}
K_0,\widetilde{K}\subset(1+O(\varepsilon))B^n.
\end{equation}
 It follows from (\ref{e14cond}) and (\ref{K0K1}) that
\begin{equation}
\label{hK0}
1\leq h_{K_0}(\tilde{u})=h_{\widetilde{K}}(\tilde{u})\leq 1+O(\varepsilon).
\end{equation}
For any $s\in(0,1)$, let $r(s)$ and $\tilde{r}(s)$
be the radii of
 $K\cap H(u,s)$ and  $\widetilde{K}\cap H(\tilde{u},s)$, respectively.
We claim that
\begin{equation}
\label{rtilder}
\tilde{r}(s)=r(s)+O(\varepsilon^{14})
\mbox{ \  if $s\leq 1-4\varepsilon^2$}.
\end{equation}
For a fixed $s\in(0,1-4\varepsilon^2]$,
let $s_1<s_2$ such that
$$
[s_1,s_2]u=\pi_{\R u}\left[K_0\cap H(\tilde{u},s)\right].
$$
Since $K_0\subset B^n+\R u$, it follows  that
\begin{equation}
\label{s1s2}
s-2\varepsilon^{16}<s_1<s_2<s+2\varepsilon^{16}.
\end{equation}
Since $1-s\geq 4\varepsilon^2$ and $u\in K_0$, we deduce that
$$
\|z-s\tilde{u}\|=r(s)+O(\varepsilon^{14})
\mbox{ \ for any $z\in \partial K_0\cap H(\tilde{u},s)$,}
$$
which in turn yields (\ref{rtilder}).

Now let $E$ be any $o$-symmetric ellipsoid having $\R u$ as an
axis of rotation. For some orthogonal linear
transform $\Phi_*$ that maps $\tilde{u}$ into $u$,
we consider the  $o$-symmetric ellipsoid
$E_*=\Phi_*^{-1}\widetilde{\Phi}^{-1}E$
 having again $\R u$ as an
axis of rotation. We know that
there exists $x_*$ such that
$x_*+4\varepsilon^2B^n\subset K_0\Delta E_*$
and $ x_*\cdot u\leq 1-4\varepsilon^2$.
It follows from (\ref{rtilder}) that
for $\tilde{x}=\Phi_*x_*$ and $\widetilde{E}=\Phi_*E_*$,
we have $\tilde{x}+3\varepsilon^2\subset \widetilde{K}\Delta\widetilde{E}$
and $ \tilde{x}\cdot \tilde{u}\leq 1-4\varepsilon^2$.
We conclude using (\ref{hK0}) and (\ref{rtilder}) that
 $x+2\varepsilon^2\,B^n\subset E\Delta K'$
and $| x\cdot u |\leq 1-2\varepsilon^2$
 for $x=\widetilde{\Phi}\tilde{x}$, verifying (i).

To prove (ii), let
$$
\varepsilon^{32}/4<p<4\varepsilon^{32}.
$$
If $tu\in H(\tilde{u},h_{K_0}(\tilde{u})-p)\cap {\rm int}\,K_0$ for $t>0$,
then
$H(\tilde{u},h_{K_0}(\tilde{u})-p)$ cuts of a cap
of depth at most $p/\varepsilon^{16}<4\varepsilon^{16}$
from $H(u,t)\cap K_0$, and hence
$H(\tilde{u},h_{K_0}(\tilde{u})-p)\cap K_0\cap H(u,t)$
is an  $(n-2)$-ball of radius at most $O(\varepsilon^{8})$.
As $K_0\subset 2B^n$, we deduce that
\begin{align*}
\mathcal{H}\left( H(\tilde{u},h_{K_0}(\tilde{u})-p)\cap \widetilde{K}\right)&=
\mathcal{H}\left( H(\tilde{u},h_{K_0}(\tilde{u})-p)\cap K_0\right)\\
&= O(\varepsilon^{8(n-2)})=O(\varepsilon^{4(n-1)}),
\end{align*}
thus for $\tilde{v}\in S^{n-1}\cap\tilde{u}^\bot$, we have
$$
(h_{K_0}(\tilde{u})-p)\tilde{u}+\gamma\varepsilon^4 \tilde{v}
\not\in \widetilde{K},
$$
where $\gamma>0$ depends on $n$.
We conclude using again (\ref{hK0}) and (\ref{rtilder}) that
$$
(1-q)u+2\gamma\varepsilon^4 v\not\in K'
\mbox{ \ for any $q\in(\varepsilon^{32}/2,2\varepsilon^{32})$},
$$
which in turn yields (ii).
\theend

Finally, we are in a position to prove Theorem~\ref{equatorpole}.\\

\noindent{\it Proof of Theorem~\ref{equatorpole}: } We assume that $\delta_0$ (and hence $\delta$, as well)
is small enough to make the estimates below work.
 We write $\gamma_1,\gamma_2,\ldots$ to denote
positive constants depending only on $n$.

 According to  Lemma~\ref{rounding} and Proposition~\ref{polefinal}, there exists
 a convex body $K_1$ spinning around $u$  and obtained from
$K$ by a combination of Steiner symmetrizations, linear transformations
and taking limits, such that  for some $\eta\in(\delta^3,\delta]$,
we have  $\delta_{\rm BM}(K_1,B^n)\leq\eta$, and
\begin{description}
\item{(a)}
for any $o$-symmetric ellipsoid
$E$ with axial rotational symmetry around $\R u$, one finds
a  ball $x+2\eta^2\,B^n\subset E\Delta K_1$
where $| x\cdot u |\leq 1-2\eta^2$;
\item{(b)} $(1-\eta^{32})u+\eta^3 v\not\in K_1$.
\end{description}
In particular,
$$
\delta_{BM}(K_1,B^n)\geq \gamma_1\eta^2.
$$

 If
$$
 \delta^3+(1-\delta^7)\,v\not\in K_1,
$$
then we simply take $\varepsilon=\eta$ and $K'=K_1$. If
$$
 \delta^3+(1-\delta^7)\,v\in K_1,
$$
then let $K_2$ be the Schwarz rounding of $K_1$
around $\R v$, and hence
$\delta_{\rm BM}(K_2,B^n)\geq \gamma_2\eta^7$
by  Lemma~\ref{equator}. For $\varepsilon=\delta_{\rm BM}(K_2,B^n)$,
we have
$$
\delta^{24} \leq\delta_{\rm BM}(K_2,B^n)=\varepsilon\leq\delta.
$$
Since $K_1\subset (1+\gamma_2 \varepsilon) B^n$ and
$K_1$ spins around $u$, if $t\in(0,\varepsilon)$, then
\begin{eqnarray*}
\mathcal{H}\left( K_1\cap H(v,1-t)\right)&\leq &
\gamma_3\varepsilon^{1/2} \mathcal{H}\left( B^n\cap u^\bot\cap H(v,1-t)\right)\leq
\gamma_4 \varepsilon^{1/2} t^{\frac{n-2}2},\\
\mathcal{H}\left(H(v,t)\cap K_1\right)&\leq & (1-\gamma_5 t^2)\mathcal{H}\left(H(v,0)\cap K_1\right).
\end{eqnarray*}
Using that $\frac{n-2}{2(n-1)}\geq\frac14$ for $n\geq 3$, we have
\begin{eqnarray*}
\gamma_6  \varepsilon^{\frac1{2(n-1)}} t^{1/4}u+(1-t)\,v&\not\in &K_2,\\
(1-\gamma_7 t^2)u+t\,v &\not\in & K_2.
\end{eqnarray*}
We transform $K_2$ into a convex body $K'$ spinning
around $u$ by a linear map, which sends $v$ into $u$,
and $v^\bot\cap K_2$ into $u^\bot\cap B^n$. We deduce that if  $t\in(0,\varepsilon/2)$, then
\begin{eqnarray}
\label{K1sec1}
(1-t)\,u+\gamma_8  \varepsilon^{\frac1{2(n-1)}} t^{1/4}\,v&\not\in &K',\\
\label{K1sec2}
t\,u+(1-\gamma_9 t^2)v &\not\in & K'.
\end{eqnarray}

In (\ref{K1sec2}), we choose $t$ such that $\varepsilon^7=\gamma_9 t^2$,
and hence
$$
\varepsilon^3\,u+(1-\varepsilon^7)\,v \not\in K'.
$$
We also deduce by substituting $t>0$
with  $\varepsilon^3=\gamma_8  \varepsilon^{\frac1{2(n-1)}} t^{1/4}$ in  (\ref{K1sec1})
that
$$
(1-\varepsilon^{32})u+\varepsilon^3 v\not\in K'.
$$
Finally suppose that for some $o$-symmetric ellipsoid
$E$ with axial rotational symmetry around $\R u$, there is no
 ball of the form $x+2\varepsilon^2\,B^n$ in
$E\Delta K'$, where $| x\cdot u |\leq 1-2\varepsilon^2$.
It follows from Corollary~\ref{pole0} that
\begin{equation}
\label{K1sec3}
(1-14\varepsilon^2)\,u+\gamma_{10}\varepsilon^{1/2}\,v \not\in K'.
\end{equation}
If $\delta_0$ is small enough, then substituting $t=14\varepsilon^2$ in
 (\ref{K1sec1}) contradicts  (\ref{K1sec3}).
Therefore $K'$ satisfies all requirements of
Theorem~\ref{equatorpole}.
\theend

\noindent {\bf Acknowledgement: } I am grateful for fruitful discussions
with Gerg\H{o} Ambrus, Erwin Lutwak, Rolf Schneider, Gaoyong Zhang during a stay at NYU-Poly. Very special thanks are due to the anonymus referees, whose remarks gave a whole new direction to the paper, and helped it to be much more clear and readable.

\noindent K.J.  B\"or\"oczky\\ 
Alfr\'ed R\'enyi Institute of Mathematics,
	 Hungarian Academy of Sciences,
	 1053 Budapest, Re\'altanoda u. 13-15.
	 HUNGARY. and

\noindent Central European University,
1051 Budapest, N\'ador u. 9, HUNGARY

\end{document}